\documentclass[3p,times]{elsarticle}

\usepackage{ecrc}
\usepackage{footnote}

\volume{00}

\firstpage{1}

\journalname{Computer Physics Communications}

\runauth{}

\jid{procs}

\jnltitlelogo{Computer Physics Communications}

\CopyrightLine{2011}{Published by Elsevier Ltd.}

\usepackage{amssymb}
\usepackage{amssymb,amsmath}

\usepackage[figuresright]{rotating}
\usepackage{subfigure}

\newcommand{\xbf}{\mathbf{x}}
\newcommand{\Rbb}{\mathbb{R}}

\usepackage{epsfig}

\begin{document}

\begin{frontmatter}





\title{Computational methods for the dynamics of the nonlinear Schr\"odinger/Gross-Pitaevskii equations}


\author[label1,label4]{Xavier Antoine}
\ead{xavier.antoine@univ-lorraine.fr}
\author[label2,label5]{Weizhu Bao\footnote{Corresponding author}}
\ead{bao@math.nus.edu.sg}
\ead[url]{http://www.math.nus.edu.sg/\~{}bao/}
\author[label3]{Christophe Besse}
\ead{christophe.besse@math.univ-lille1.fr}

\address[label1]{Universit\'e de Lorraine, Institut Elie Cartan de
Lorraine, UMR 7502, Vandoeuvre-l\`es-Nancy, F-54506, France}

\address[label4]{CNRS, Institut Elie Cartan de Lorraine, UMR 7502, Vandoeuvre-l\`es-Nancy, F-54506, France}

\address[label2]{Department of Mathematics,
National University of Singapore, Singapore 119076, Singapore}
\address[label5]{Center for Computational Science and Engineering,
National University of Singapore, Singapore 119076, Singapore}

\address[label3]{Laboratoire Paul Painlev\'e, Universit\'e Lille Nord de France,
  CNRS UMR 8524, INRIA SIMPAF Team,
  Universit\'e Lille 1 Sciences et Technologies,
  Cit\'e Scientifique, 59655 Villeneuve d'Ascq Cedex, France}

\begin{abstract}
In this paper, we begin with the nonlinear Schr\"{o}dinger/Gross-Pitaevskii  equation (NLSE/GPE) for modeling
Bose-Einstein condensation (BEC) and nonlinear optics as well as other applications,
and discuss their dynamical properties ranging from
time reversible, time transverse invariant, mass and energy conservation, dispersion relation
to soliton solutions. Then, we review and compare different numerical methods for solving the NLSE/GPE
including finite difference time domain methods and time-splitting spectral method, and discuss
different absorbing  boundary conditions. In addition, these numerical methods are extended
to the NLSE/GPE  with damping terms and/or an angular momentum rotation term as well
as coupled NLSEs/GPEs. Finally, applications to simulate a quantized vortex lattice
dynamics in a rotating BEC are reported.
\end{abstract}


\end{frontmatter}

\section{Introduction} \setcounter{equation}{0}

The nonlinear Schr\"{o}dinger equation (NLSE) is a partial differential equation (PDE)
that can be met in many different areas of physics and chemistry as well as engineering
\cite{ADP,AS,DP,Dirac,PitaevskiiStringari,Sch,SulemSulem}.
The most well-known and important derivation of the NLSE is from the mean-field approximation
of many-body problems in quantum physics and chemistry \cite{PitaevskiiStringari},
especially for the study of Bose-Einstein
condensation (BEC) \cite{Anderson,PitaevskiiStringari}, where it is also called  the
Gross-Pitaevskii equation (GPE) \cite{BaoCai313,Gross,Pitaevskii,PitaevskiiStringari}.
Another important application of the NLSE is for laser beam propagation in nonlinear and/or quantum
optics and  there it is also known as parabolic/paraxial approximation of the
Helmholtz or time-independent Maxwell equations \cite{ADP,AS,DP,SulemSulem}. Other important applications
include long range wave propagation in underwater acoustics \cite{N,SulemSulem}, plasma and particle physics \cite{SulemSulem},
semiconductor industry \cite{Mark,MRS2002}, materials simulation based on the first principle \cite{Eng},
superfluids \cite{Aft,Bar},
molecular systems in biology \cite{DA}.

In this paper, we  consider the following time-dependent
NLSE \cite{BaoCai313,SulemSulem}
\begin{equation}\label{eqSchrodinger}
i \varepsilon \partial_t \psi(t,\xbf) = -\frac{\varepsilon^2}{2} \nabla^2 \psi(t,\xbf)
+V(\xbf) \psi(t,\xbf) +f(|\psi(t,\xbf)|^2) \psi(t,\xbf), \qquad \xbf \in \mathbb{R}^d,\quad t>0,
\end{equation}
with initial data
\begin{equation} \label{init0}
\psi(t=0,\xbf)=\psi_0(\xbf),\qquad \xbf \in \mathbb{R}^d,
\end{equation}
where $i=\sqrt{-1}$ is the complex unit, $t$ is the time variable,
$\xbf\in\mathbb{R}^d$ is the spatial variable with $d=1, 2, 3$,
$\psi:=\psi(t,\xbf)$ is the complex-valued wave function or order parameter,
$\nabla^2=\Delta$ is the usual Laplace operator and  $\psi_{0}:=\psi_{0}(\xbf)$ is a given
complex-valued initial data. Moreover, $0<\varepsilon\leq 1$ is a dimensionless parameter. In most physics literatures,
it is taken equal to $\varepsilon=1$; however in the semi-classical regime, we have $0<\varepsilon \ll 1$,
$\varepsilon$ being  also called  the ''scaled'' Planck constant.  Function $V:=V(\xbf)$ is a given
real-valued external potential and its specific form depends on different applications
and could also sometimes be time dependent \cite{Bao01,BaoCai313,PitaevskiiStringari,SulemSulem}. For example, in BEC,
it is usually chosen as either a harmonic confining potential, i.e. $V(\xbf)= \frac{|\xbf|^{2}}{2}$
\cite{Bao01,BaoCai313,PitaevskiiStringari} or
an optical lattice potential, i.e.  $V(\xbf)= A_{1}\cos(L_{1} x)+A_{2}\cos(L_{2} y)+A_{3}\cos(L_{3} z)$
in three dimensions (3D) with $A_1$, $A_2$, $A_3$, $L_1$, $L_2$ and $L_3$ constants \cite{Bao01,BaoCai313,Huang,Huang1,PitaevskiiStringari},
or a stochastic potential for producing speckle patterns \cite{MinLRB,PitaevskiiStringari};  in nonlinear optics,
it might be chosen as an attractive potential, i.e. $V(\xbf)= - \frac{|\xbf|^{2}}{2}$ \cite{ADP,DP}.
The nonlinearity $f:=f(\rho)$ is a real-valued smooth function depending on the density
$\rho:=|\psi|^{2}\in[0,\infty)$ and its specific form depends on different applications \cite{ADP,AS,Bao01,BaoCai313,PitaevskiiStringari,SulemSulem}.
In fact, when $f(\rho)\equiv\lambda$, with $\lambda$ a constant, then the NLSE (\ref{eqSchrodinger})
collapses to the standard (linear)
Schr\"{o}dinger equation \cite{Dirac,Sch}. The most popular and important nonlinearity
is the cubic nonlinearity \cite{ADP,AS,BaoCai313,PitaevskiiStringari,SulemSulem}
\begin{equation}
\label{nls_cubic}
f(\rho)= \beta \rho, \qquad 0\le \rho<\infty,
\end{equation}
where $\beta$ (positive for repulsive or defocusing interaction and negative for attractive or
focusing interaction) is a given dimensionless constant describing
the strength of interaction. Other nonlinearities used in nonlinear optics
include the cubic-quintic nonlinearity $f(\rho)=\beta_{1} \rho + \beta_{2}\rho^2$ \cite{ADP,DP,N,SulemSulem}
and the saturation of the intensity
nonlinearity  $f(\rho)=\frac{\beta_{0}\rho}{1+c_{0}\rho}$ with
$\beta_0$, $\beta_1$, $\beta_2$ and $c_0$ given constants \cite{ADP,DP,N,SulemSulem}.
In some applications, nonlocal nonlinearities can also be considered, and in this case
$f(\rho)=U*\rho$ which is a convolution, with $U:=U(\xbf)$ the kernel \cite{BaoCai313,BaoCaiWang,BMS,BEGMY,Cai,Xiong}. Specifically,
for the Hartree potential case, $U(\xbf)=\frac{1}{4\pi |\xbf|}$ is
the Green's function of the Laplace operator in 3D, and then,
the NLSE (\ref{eqSchrodinger}) can be re-written as the Schr\"{o}dinger-Poisson system
in 3D \cite{BEGMY,BaoJMZ}. Finally, we remark here that the NLSE (\ref{eqSchrodinger}) is also
called  the Gross-Pitaevskii Equation (GPE) when the nonlinearity is chosen as the cubic
nonlinearity as in (\ref{nls_cubic}), especially in BEC \cite{Bao01,BaoCai313,PitaevskiiStringari}.
Thus GPE is a special version of NLSE \cite{Bao01,BaoCai313,Gross,Pitaevskii,PitaevskiiStringari}.

There are many important dynamical properties of the solution $\psi$ to the NLSE (\ref{eqSchrodinger}).
Among them, here we mention several important
ones that we will use  to justify different numerical methods
on whether they are still valid at the discrete level after the NLSE (\ref{eqSchrodinger}) is
discretized by a numerical method. In fact, the NLSE (\ref{eqSchrodinger}) is a dispersive
PDE and it is {\sl time reversible or symmetric}, i.e.
it is unchanged under the change of variable in time as $t \rightarrow -t$ and taken conjugate in the equation.
Another important property is {\sl time transverse or gauge invariant}, i.e.
if $V \rightarrow V+\alpha$, with $\alpha$ a given real constant,
then the solution $\psi \rightarrow \psi e^{-i\alpha t/\varepsilon}$ which immediately implies that
the density $\rho=|\psi|^2$ is unchanged. The NLSE (\ref{eqSchrodinger}) conserves many quantities.
Among them, the {\sl mass} (or wave energy in nonlinear optics) and {\sl energy} (or Hamiltonian in nonlinear optics)
{\sl conservation} are given as \cite{BaoCai313,Cazenave,SulemSulem}
\begin{eqnarray}\label{massConserveC}
&&N(t):=\|\psi(t,\cdot)\|^2=\int_{\mathbb{R}^{d}} |\psi(t,\xbf)|^{2} d \mathbf{x} \equiv \int_{\mathbb{R}^{d}} |\psi(0,\xbf)|^{2} d \mathbf{x} = \int_{\mathbb{R}^{d}} |\psi_0(\xbf)|^{2} d \mathbf{x}:=N(0), \qquad
 t \ge 0, \\
&&E(t):=\int_{\mathbb{R}^{d}} \left[\frac{\varepsilon^2}{2} |\nabla \psi(t,\xbf) |^{2} + V(\xbf) | \psi(t,\xbf)|^{2} +F(|\psi(t,\xbf)|^2)\right] d \mathbf{x} \equiv E(0), \qquad  t \ge 0,
\label{energyConserveC}
\end{eqnarray}
respectively, with
\begin{equation}\label{Frho}
F(\rho):=\int_0^{\rho} f(s)\,ds, \qquad 0\le \rho<+\infty.
\end{equation}

If there is no external potential in the NLSE (\ref{eqSchrodinger}), i.e.  $V(\xbf)\equiv0$,
then the momentum and angular momentum are also conserved \cite{BaoCai313,Cazenave,SulemSulem}.
In addition, the NLSE (\ref{eqSchrodinger}) admits the plane wave solution as
$\psi(t,\xbf) = A e^{i(\mathbf{k} \cdot \xbf - \omega t)}$, where
the time frequency $\omega$, amplitude $A$ and spatial wave number $\mathbf{k}$ satisfy the following
{\sl dispersion relation} \cite{BaoCai313,Cazenave,Ignat,SulemSulem}
\begin{equation}\label{dispre}
\omega = \frac{\varepsilon|\mathbf{k}|^{2}}{2}+\frac{1}{\varepsilon}f(|A|^{2}).
\end{equation}
In this case, if in 1D with $\varepsilon=1$ and the nonlinearity is chosen as the
focusing cubic nonlinearity (\ref{nls_cubic}) with $\beta <0$, it also
admits the well-known bright soliton solution as \cite{ADP,AS,BaoTangXu,DP,N,SulemSulem}
\begin{equation}
\label{bright_sol}
\psi_{B}(t,x)=\frac{A}{\sqrt{-\beta}} \textrm{sech}(A(x-vt-x_{0})) e^{i(vx -\frac{1}{2}(v^{2}-A^{2})t + \theta_{0})},
\qquad x \in \Rbb, \quad t \geq 0,
\end{equation}
where $\frac{A}{\sqrt{-\beta}}$ is the amplitude of the soliton with $A$ a real
constant, $v$ is the velocity of the soliton, $x_{0}$
and $\theta_0$ are the initial shifts in space and phase, respectively.
Since the soliton solution is exponentially decaying
 for $|x| \rightarrow + \infty$, then the mass and energy are well defined and given by: $N(\psi_{B})=-\frac{2A}{\beta}$
 and $E({\psi_{B}})=\frac{Av^{2}}{-\beta}+\frac{A^{3}}{-3\beta}$.
For other more dynamical properties of the NLSE (\ref{eqSchrodinger}) such as
dark (black and grey) solitons in 1D under defocusing cubic nonlinearity,
well-posedness and finite time blow-up, we refer to
\cite{BaoCai313,Cazenave,Fibich,FibichP,SulemSulem} and references therein.

For studying numerically the dynamics of the NLSE (\ref{eqSchrodinger}),
different numerical methods have been proposed in the literatures
\cite{ADK91,Review08,Bao01,BaoCai313,BaoJakschP,BaoJinP,BaoJinP2,Besse04,Chan1,Chan2,
Tosi2,Cerim,sanzserna,Guo,Hardin,Jin,MaPiPo,Path,Taha1,WH86} and references therein.
The main aim of this paper is to review different numerical methods
for solving the NLSE (\ref{eqSchrodinger}) numerically, compare them
in terms of keeping different dynamical properties at the discretized level,
stability and accuracy, discuss
different absorbing  boundary conditions for truncating
the NLSE (\ref{eqSchrodinger}) onto a bounded computational domain,
 and extend them for solving
NLSE with damping and/or angular momentum terms as well as coupled NLSEs.

The paper is organized as follows. In Section \ref{s2}, we review several
popular numerical methods in the literatures for discretizing NLSE/GPE under simple boundary conditions
 and compare their advantages and disadvantages, discuss different absorbing
boundary conditions for NLSE/GPE, and review some robust and efficient numerical methods
for NLSE/GPE in the semiclassical regime. Extensions to NLSE/GPE with damping and/or angular
rotating terms and coupled NLSEs/GPEs are presented in Sections \ref{s3} and \ref{s4},
respectively. Numerical comparison of different numerical methods and some applications
are reported in Section \ref{s5}. Finally, we  draw some
conclusions and discuss future perspectives in Section \ref{s6}.

\section{Numerical methods for NLSE/GPE}\label{s2}
\setcounter{equation}{0}

\subsection{Some popular numerical methods}
Here we present several popular numerical methods for discretizing the NLSE/GPE (\ref{eqSchrodinger})
 and discuss/compare their properties including stability, accuracy, computational cost and how
much properties of NLSE/GPE that these numerical methods can keep at the discretized level.
For simplicity of notation, we shall introduce these methods  in one
space dimension (1D), i.e. $d=1$ in (\ref{eqSchrodinger}).  Generalizations to $d>1$ are straightforward for
tensor product grids and the results remain valid with modifications.
For $d=1$, the NLSE/GPE (\ref{eqSchrodinger}) truncated on a bounded interval $(a,b)$ with homogeneous
Dirichlet boundary condition ($|a|$ and $b$ large enough so that the error due to the truncation
is negligible) becomes
\begin{eqnarray}\label{NLSE1}
&&i \varepsilon \partial_t \psi(t,x) = -\frac{\varepsilon^2}{2} \partial_{xx} \psi(t,x)
+V(x) \psi(t,x)  +f(|\psi(t,x) |^2) \psi(t,x) , \qquad a<x<b,\quad t>0,\\
\label{NLSE2}
&&\psi(t,a)=\psi(t,b)=0, \qquad t\ge0,\\
 \label{NLSE3}
&&\psi(t=0,x)=\psi_0(x),\qquad a\le x\le b.
\end{eqnarray}
In this case, the {\sl mass} and
{\sl energy conservation} collapse to the following
\begin{eqnarray}\label{massConserveC1d}
&&N(t):=\|\psi(t,\cdot)\|^2=\int_a^b |\psi(t,x)|^{2} dx \equiv \int_a^b |\psi(0,x)|^{2} dx
= \int_a^b |\psi_0(x)|^{2} dx:=N(0), \qquad t \ge 0, \\
&&E(t):=\int_a^b \left[\frac{\varepsilon^2}{2} |\partial_x \psi(t,x) |^{2} + V(x) | \psi(t,x)|^{2} +F(|\psi(t,x)|^2)\right] dx \equiv E(0), \qquad t \ge 0.
\label{energyConserveC1d}
\end{eqnarray}

Choose a time step $\tau:=\Delta t>0$ and denote the different times as $t_n:=n\tau$ for $n=0,1,2,\ldots$;
choose the mesh size $h :=\Delta x=\frac{b-a}{J}$, with $J$ an even positive integer, and denote the
grid points by $x_{j}:=a+jh$, for $j=0,1,\ldots,J$. Let $\psi_j^n$ be the numerical approximation
of $\psi(t_n,x_j)$, for $j=0,1,\ldots,J$ and $n=0,1,2,\ldots$; let $\boldsymbol{\psi}^n$ be the solution vector
at time $t=t_n$ with components $(\psi_j^n)_{0\le j\le J}$. In addition,
for any complex-valued  vector $\boldsymbol{\phi} = (\phi_j)_{0\le j\le J}$, we define the following
finite difference operators
\begin{equation}
(D^+_x \boldsymbol{\phi})_j=\frac{\phi_{j+1}-\phi_j}{h}, \quad 0\le j\le J-1;\qquad \qquad
 (D^2_x \boldsymbol{\phi})_j=\frac{\phi_{j+1}-2\phi_j+\phi_{j-1}}{h^2}, \qquad 1\le j\le J-1.
\end{equation}

Then the  {\sl Crank-Nicolson finite difference} (CNFD) method --- in which one applies the second-order centered
difference scheme for spatial discretization and the Crank-Nicolson scheme  for time discretization ---
for discretizing (\ref{NLSE1}) reads as \cite{BaoCai313,BaoCai,Chang1,Chang,Glassey,Guo,Zhu}
\begin{eqnarray}
\label{fulldiscreteCNFD}
i\varepsilon\frac{\psi^{n+1}_j-\psi^{n}_j}{\tau}=- \frac{\varepsilon^2}{4}\left[\left(D^2_x\boldsymbol{\psi}^{n+1}\right)_j+
\left(D^2_x\boldsymbol{\psi}^n\right)_j\right]+\left[V(x_j)+
G\left(|\psi_j^{n+1}|^2,|\psi_j^{n}|^2\right)\right]\frac{\psi_j^{n+1}+\psi_j^n}{2},  \quad 1\le j\le J-1,\ \
n\ge0,
\end{eqnarray}
where
\[
{G}(\rho_{_1},\rho_{_2})=\frac{F(\rho_{_1})-F(\rho_{_2})}{\rho_{_1}-\rho_{_2}}:=\int_0^1 f(\theta \rho_{_1}+
(1-\theta)\rho_{_2})\,d\theta, \qquad 0\le \rho_{_1},\rho_{_2}<\infty.
\]
The boundary and initial conditions (\ref{NLSE2})-(\ref{NLSE3}) are discretized as
\begin{equation}\label{FDmatrix}
\psi_0^{n+1}=\psi_J^{n+1}=0, \qquad n\ge0; \qquad \qquad \psi_j^{0}=\psi_0(x_j), \qquad 0\le j\le J.
\end{equation}
The above CNFD method (\ref{fulldiscreteCNFD}) is {\sl time reversible or symmetric}, i.e.
it is unchanged if $\tau\to -\tau$ and $n+1\leftrightarrow n$; and its memory cost is $\mathcal{O}(J^d)$
in $d$-dimensions with $J$ unknowns in each direction.
It conserves the {\sl mass}
(\ref{massConserveC1d}) and {\sl energy} (\ref{energyConserveC1d}) in the discretized level \cite{BaoCai313,BaoCai,BaoCai2,Chang,Glassey}, i.e.
\begin{eqnarray} \label{massdis}
&&N^{n}:=h \sum_{j=1}^{J-1} |\psi_j^n|^2 \equiv
h \sum_{j=1}^{J-1} |\psi_j^0|^2 =h \sum_{j=1}^{J-1} |\psi_0(x_j)|^2:=N^0, \qquad n\ge0,\\
\label{energydis}
&&E^{n}:=  h\sum_{j=0}^{J-1}\left[\frac{\varepsilon^2}{2}\left|\left(D_x^+ \boldsymbol{\psi}^{n}\right)_j\right|^2
+V(x_j)|\psi_j^n|^2+F(|\psi_j^n|^2)\right]\equiv E^0, \qquad n\ge0,
\end{eqnarray}
which immediately implies that the CNFD method is unconditionally stable. In addition,
it can be proven rigorously in mathematics that the CNFD method is second-order accurate in both
time and space for any fixed $\varepsilon=\varepsilon_0=\mathcal{O}(1)$
\cite{BaoCai313,BaoCai,BaoCai2,Chang,Glassey}. However, it cannot preserve the {\sl time transverse invariant} and
{\sl dispersion relation} properties of the NLSE/GPE (\ref{NLSE1}) at the discretized level
\cite{BaoCai313,BaoCai,BaoCai2}. Furthermore, it is an implicit scheme and
its practical implementation is a little tedious.
In fact, at each time step,
one needs to solve a coupled fully nonlinear system
which can be solved by a fixed point or a modified Newton-Raphson  iterative
method \cite{A93,AD91,ADK91,BaoCai313,BaoCai,BaoCai2} and thus it might be very time consuming. In general, the computational cost
per time step is  much larger than $\mathcal{O}(J^d)$, especially in 2D and 3D.
In fact, if the nonlinear system in (\ref{fulldiscreteCNFD}) is not solved very
accurately, e.g. up to machine precision, the numerical solution computed from (\ref{fulldiscreteCNFD})
does not conserve the energy in (\ref{energydis}) {\sl exactly} \cite{BaoCai2}.

Due to the high computational cost of the CNFD method, in the literatures
it is motivated for considering some
semi-implicit methods in which a linear system is to be solved at every time step. Thus, the
computational cost is significantly reduced, especially in 2D and 3D.  One of  these
numerical methods is the {\sl relaxation finite difference} (ReFD) method --- in which one
gives a name to the nonlinear term and makes a second-order approximation at time $t_n$ ---
for the NLSE (\ref{NLSE1})
\cite{Besse04}
\begin{equation}
\left \{
\begin{array}{l}
\displaystyle
\frac{1}{2} \left({u}^{n+1/2}_j+{u}^{n-1/2}_j\right)
 = {f}(|{\psi}^{n}_j|^{2})\, ,\qquad \qquad \qquad \qquad \qquad \qquad 1\le j\le J-1,\qquad
n\ge0,\\[2mm]
\displaystyle
i\varepsilon  \frac{{\psi}^{n+1}_j-{\psi}^{n}_j}{\tau}
  =-\frac{\varepsilon^2}{4}\left[\left(D^2_x\boldsymbol{\psi}^{n+1}\right)_j+
\left(D^2_x\boldsymbol{\psi}^n\right)_j\right] + \frac{1}{2}\left[V(x_j) + {u}^{n+1/2}_j \right]
\left(\psi_j^{n+1}+\psi_j^n\right),
\label{Relaxsemidiscrete}
\end{array}
\right .
\end{equation}
where $u^{-1/2}_j=\psi^{0}_j=\psi_0(x_j)$ for $0\le j\le J$. The boundary and initial conditions (\ref{NLSE2})-(\ref{NLSE3}) are discretized as (\ref{FDmatrix}).
Similarly, this ReFD method is {\sl time reversible or symmetric}, its memory cost is $\mathcal{O}(J^d)$
in $d$-dimensions with $J$ unknowns in each direction, and it is easier to be implemented
than that of CNFD.
It is an implicit scheme where, at each time step, only a linear system needs to be solved for example
by the Thomas algorithm  at the cost of $\mathcal{O}(J)$ in 1D and
by some iterative methods such as conjugate gradient (CG) or multigrid (MG) method \cite{Besse04,Chang3}
in 2D and 3D
at the cost, in general, larger than $\mathcal{O}(J^d)$.
Thus it is computationally much cheaper than that of the CNFD method.
It conserves the {\sl mass}
(\ref{massConserveC}) at the discretized level as (\ref{massdis}) \cite{Besse04}, and thus
it is unconditionally stable. In addition,
it can be mathematically proven rigorously that the ReFD method is second-order accurate in both
time and space for any fixed $\varepsilon=\varepsilon_0=\mathcal{O}(1)$ \cite{Besse04}. Furthermore, for the NLSE with the cubic nonlinearity (\ref{nls_cubic}),
this method also conserves the following discrete energy defined as
\begin{equation}
\tilde{E}^{n}:=h\sum_{j=0}^{J-1}\frac{\varepsilon^2}{2}\left|\left(D_x^+ \boldsymbol{\psi}^{n}\right)_j\right|^2
+h\sum_{j=1}^{J-1}\left[V(x_j)|\psi_j^n|^2+\frac{\beta}{2}u^{n+1/2}_ju^{n-1/2}_j\right]\equiv \tilde{E}^0,
\qquad n \geq 0.
\label{energyDefDiscreteRelax}
\end{equation}
Of course, it cannot preserve the {\sl time transverse invariant} and
{\sl dispersion relation} properties of the NLSE/GPE (\ref{NLSE1}) in the discretized level
\cite{Besse04}. In addition, for general nonlinearity other than the cubic nonlinearity, no discrete energy
conservation has been proven for the ReFD method yet. Another popular method is the
following {\sl semi-implicit finite difference} (SIFD) method --- in which one uses the leap-frog scheme
for the nonlinear term  in time discretization ---
for the NLSE (\ref{NLSE1}) \cite{BaoCai,BaoCai313,BaoCai2}
\begin{eqnarray}
\label{SIFD}
i\varepsilon\frac{\psi^{n+1}_j-\psi^{n-1}_j}{2\tau}=- \frac{\varepsilon^2}{4}\left[\left(D^2_x\boldsymbol{\psi}^{n+1}\right)_j+
\left(D^2_x\boldsymbol{\psi}^{n-1}\right)_j\right]+\left[V(x_j)+
f(|\psi_j^{n}|^2)\right]\psi_j^n,  \qquad 1\le j\le J-1,\quad n\ge1;
\end{eqnarray}
or
\begin{eqnarray}
\label{SIFD987}
i\varepsilon\frac{\psi^{n+1}_j-\psi^{n-1}_j}{2\tau}=- \frac{\varepsilon^2}{4}\left[\left(D^2_x\boldsymbol{\psi}^{n+1}\right)_j+
\left(D^2_x\boldsymbol{\psi}^{n-1}\right)_j\right]+V(x_j)\frac{\psi_j^{n+1}+\psi_j^n}{2}+
f(|\psi_j^{n}|^2)\psi_j^n,  \quad 1\le j\le J-1,\ n\ge1;
\end{eqnarray}
The boundary and initial conditions (\ref{NLSE2})-(\ref{NLSE3}) are
discretized as (\ref{FDmatrix}) and the first step can be computed as \cite{BaoCai313,BaoCai,BaoCai2}
\begin{equation}\label{fstep87}
\psi_j^1=\psi_j^0-\frac{i\tau}{\varepsilon}\left[-\frac{\varepsilon^2}{2}
\left(D^2_x\boldsymbol{\psi}^{0}\right)_j+V(x_j)\psi_j^0+f(|\psi_j^0|^2)\psi_j^0\right],
\qquad 1\le j\le J-1.
\end{equation}
Again, this SIFD method is {\sl time reversible or symmetric}, i.e.
it is unchanged if $\tau\to -\tau$ and $n+1\leftrightarrow n-1$; its memory cost is $\mathcal{O}(J^d)$
in $d$-dimensions considering $J$ unknowns in each direction. It is also much easier to be implemented
than that of CNFD and ReFD, especially in 2D and 3D.
It is an implicit scheme where, at each time step, only a linear system (whose coefficient matrix
is time independent) needs to be solved,
 by the Thomas algorithm  in  $\mathcal{O}(J)$  operations in 1D and
by the direct fast Poisson solver \textit{via} the discrete sine transform (DST) at the cost
of $\mathcal{O}(J^d\ln J)$ in 2D and 3D. Thus, it is significantly cheaper than the
CNFD and ReFD methods, especially in 2D/3D. In addition,
it can be proven that the SIFD method is second-order accurate  both in
time and space for any fixed $\varepsilon=\varepsilon_0=\mathcal{O}(1)$
\cite{BaoCai313,BaoCai,BaoCai2} and a fourth-order (or six-order or even eighth-order) accurate method
can be very easily constructed \textit{via} the Richardson extrapolation technique \cite{WangTC}.
Of course, it is conditionally stable under the stability condition
$\tau\lesssim \tau_n:=\frac{1}{\varepsilon \max_{0\le j\le J} |V(x_j)+f(|\psi_j^n|^2)|}$
(or $\tau\lesssim \tau_n:=\frac{1}{\varepsilon \max_{0\le j\le J} |f(|\psi_j^n|^2)|}$ for (\ref{SIFD987}))
for $n\ge0$, which is independent of the mesh size $h$.
It cannot preserve the {\sl time transverse invariant},
{\sl dispersion relation} and {\sl mass} and {\sl energy conservation}
properties of the NLSE/GPE (\ref{NLSE1}) in the discretized level
\cite{BaoCai313,BaoCai,BaoCai2}.

Another way to handle the nonlinearity in (\ref{NLSE1}) is to use the time-splitting
technique \cite{BaoCai313,BaoJakschP,BaoJinP,BaoJinP2,BaoShen,Huang,Huang1,
Thalhammer3,Hardin,Path,Taha1,Thalhammer,WH86}, i.e.
from time $t=t_n$ to time $t=t_{n+1}$, the equation (\ref{NLSE1}) is solved in two splitting steps.
One solves first
\begin{eqnarray}\label{NLSE1st}
i \varepsilon \partial_t \psi(t,x) = -\frac{\varepsilon^2}{2} \partial_{xx} \psi(t,x),
\qquad a<x<b, \qquad t>t_n,
\end{eqnarray}
with the homogeneous Dirichlet boundary condition (\ref{NLSE2})
for the time step of length $\tau$, followed by solving
\begin{eqnarray}\label{NLSE2nd}
i \varepsilon \partial_t \psi(t,x) =
V(x) \psi(t,x)  +f(|\psi(t,x) |^2) \psi(t,x) , \qquad a\le x \le b,\quad t>t_n,
\end{eqnarray}
for the same time step. Eq. (\ref{NLSE1st}) will be first discretized in space by the sine spectral method and
then integrated (in phase space) in time {\sl exactly} \cite{BaoCai313,BaoJakschP,BaoJinP,BaoJinP2,BaoZhang}.
Multiplying (\ref{NLSE2nd}) by $\overline{\psi(t,\xbf)}$
(conjugate of $\psi(t,\xbf)$) and then subtracting from its conjugate  \cite{BaoCai313,BaoJakschP,BaoJinP,BaoJinP2},
we get for $\rho(t,x):=|\psi(t,x)|^2$
\begin{equation}
\partial_t\,\rho(t,x)=0, \qquad t> t_n, \quad a\le x\le b,
\end{equation}
which immediately implies that the density $\rho$ is invariant for any fixed $x$ in the second splitting step
(\ref{NLSE2nd}), i.e.
\begin{equation}\label{rhovar}
\rho(t,x) :=|\psi(t,x)|^2\equiv|\psi(t_n,x)|^2=\rho(t_n,x), \qquad t\ge t_n, \quad
a\le x\le b.
\end{equation}
Plugging (\ref{rhovar}) into (\ref{NLSE2nd}), Eq. (\ref{NLSE2nd}) collapses to a linear ODE as
\begin{eqnarray}\label{NLSE2nd76}
i \varepsilon \partial_t \psi(t,x) =
V(x) \psi(t,x)  +f(|\psi(t_n,x) |^2) \psi(t,x) , \qquad a\le x \le b,\quad t>t_n,
\end{eqnarray}
which can be integrated {\sl analytically} as
\begin{equation}
\psi(t,x)=e^{-i\left[V(x)+f(|\psi(t_n,x) |^2)\right](t-t_n)/\varepsilon}\, \psi(t_n,x),
\qquad t\ge t_n, \quad a\le x\le b.
\end{equation}
In practical computation, from time $t=t_n$ to $t=t_{n+1}$, one often combines the splitting steps \textit{via}
the standard Strang splitting \cite{Strang} ---
which is usually referred to as time-splitting sine pseudospectral (TSSP) method
\cite{BaoCai313,BaoJakschP,BaoJinP,BaoJinP2,BaoZhang} --- as
 \begin{equation}
 \begin{array}{l}
\psi_{j}^{(1)}= e^{-i \tau [V(x_{j})+f(|\psi_{j}^{n})|^2)]/2\varepsilon }\psi_{j}^{n}\, , \\
\displaystyle \psi_{j}^{(2)}=\sum_{l=1}^{J-1} e^{-i\tau \, \varepsilon\,\mu_{l}^{2}/2} (\widehat{\psi^{(1)}} )_{l} \sin\left(\mu_l(x_j-a)\right)=\sum_{l=1}^{J-1} e^{-i\tau \, \varepsilon\,
\mu_{l}^{2}/2} (\widehat{\psi^{(1)}} )_{l} \sin\left(\frac{lj\pi}{J}\right),\\
\psi_j^{n+1}=e^{-i \tau[V(x_{j})+f(|\psi_{j}^{(2)}|^2)]/2\varepsilon}\psi_{j}^{(2)}\, ,
\qquad \qquad \qquad \qquad \qquad \qquad 1\le j\le J-1,\quad n\ge0,
 \end{array}
 \end{equation}
where $\psi_0^{n+1}=\psi_J^{n+1}=0$ for $n\ge0$, $\psi_j^0=\psi_0(x_j)$ for $0\le j\le J$, and
$(\widehat{\psi^{(1)}} )_{l}$ for $1\le l\le J-1$, the discrete sine transform coefficients of the complex-valued vector $\psi^{(1)}=(\psi^{(1)}_0,\psi^{(1)}_1,\cdots,\psi^{(1)}_J)^T$ with $\psi^{(1)}_0=\psi^{(1)}_J=0$, are defined by
 \begin{equation}
\mu_l = \frac{\pi l}{b-a}, \qquad \qquad\qquad
\widehat{\psi}^{(1)}_{l}=\frac{2}{J}\sum_{j=1}^{J-1} \psi_{j}^{(1)} \sin\left(\mu_l(x_{j}-a)\right)
=\frac{2}{J}\sum_{j=1}^{J-1} \psi_{j}^{(1)} \sin\left(\frac{lj\pi}{J}\right),
\qquad 1\le l\le J-1.
\end{equation}
Again, the above TSSP method is {\sl time reversible or symmetric}, its memory cost is $\mathcal{O}(J^d)$
in $d$-dimensions with $J$ unknowns along each direction. Its implementation is much easier
than for CNFD, ReFD and SIFD.
It is an explicit scheme since there is no need to solve a linear system and, at each time step,
the computational cost is $\mathcal{O}(J^d\ln J)$.
It conserves the {\sl mass}
(\ref{massConserveC1d}) in the discretized level as (\ref{massdis})
\cite{BaoCai313,BaoJakschP,BaoJinP,BaoJinP2,BaoZhang}, and
it is therefore unconditionally stable. In addition,
it can be rigorously proven that the TSSP method is second-order accurate in time
 and spectral order accurate in space for any fixed $\varepsilon=\varepsilon_0=\mathcal{O}(1)$
 \cite{BaoCai313,Besse,Gau,Lubich,Neuhauser,ShenWang2,Thalhammer1}. In addition, it is {\sl time transverse invariant},
i.e. when $V(x)\to V(x)+\alpha$ with $\alpha$ a constant, then $\psi_j^n\to \psi_j^n e^{-i\alpha n\tau/\varepsilon}$
which implies that the density $\rho_j^n:=|\psi_j^n|^2$ is unchanged;  it has the same
{\sl dispersion relation} as  the NLSE/GPE (\ref{NLSE1}), i.e. if $V(x)\equiv 0$ and $\psi_0(x_j) =A e^{ikx_j} $ for $0\le j\le J$, then the numerical solution from the TSSP method is $\psi_j^n=A e^{i(kx_j - \omega t_n)}$ with $\omega$, $A$ and $k$ satisfying the dispersion relation (\ref{dispre}) provided that $J\ge k$.
However, it cannot preserve  {\sl energy conservation}
properties of the NLSE/GPE (\ref{NLSE1}) in the discretized level
\cite{BaoJinP,BaoJinP2}. In fact, the numerical solution from the TSSP method actually
coincides with the exact solution of a modified PDE at each time step. This shows the existence of a modified energy \cite{Debussche,Duj0,Duj,FaouBook,Faou1} preserved by the numerical scheme that is close to the exact energy if the numerical solution is smooth. However, some resonances may destroy the energy conservation on long time evolution (see Chapter 7 of \cite{FaouBook}). Finally, for designing high-order, e.g. fourth-order accurate in time,
time-splitting spectral methods, we refer to \cite{BaoShen,BaoZhang,Thalhammer,Thalhammer1} and references therein.

 In the literatures, in some applications where the solution of the NLSE/GPE is not smooth, e.g.
with random potential \cite{And,Dubi,Lye,MinLRB}, then the following
time-splitting finite difference (TSFD) -- in which the time-splitting is
applied first and then the free Schr\"{o}dinger equation (\ref{NLSE1st}) is discretized by the CNFD method ---
is used for discretizing the NLSE/GPE \cite{BaoCai313,BaoTang,BaoTangXu,Wang0} as
\begin{equation}
 \begin{array}{l}
\psi_{j}^{(1)}= e^{-i \tau [V(x_{j})+f(|\psi_{j}^{n}|^2)]/2\varepsilon }\psi_{j}^{n}\, , \qquad 0\le j\le J,\\
\displaystyle i\varepsilon \frac{\psi_{j}^{(2)}-\psi_{j}^{(1)}}{\tau}=
-\frac{\varepsilon^2}{4}\left[\left(D^2_x\boldsymbol{\psi}^{(2)}\right)_j+
\left(D^2_x\boldsymbol{\psi}^{(1)}\right)_j\right], \qquad 1\le j\le J-1, \qquad \qquad
\psi_{0}^{(2)}=\psi_{J}^{(2)}=0,\\
\psi_j^{n+1}=e^{-i \tau[V(x_{j})+f(|\psi_{j}^{(2)}|^2)]/2\varepsilon}\psi_{j}^{(2)}\, ,
\qquad \qquad  0\le j\le J,\quad n\ge0,
 \end{array}
 \end{equation}
where $\psi_j^0=\psi_0(x_j)$ for $0\le j\le J$. Similarly,  the above TSFD method
is {\sl time reversible or symmetric}, its memory cost is $\mathcal{O}(J^d)$
in $d$-dimensions with $J$ unknowns in each direction; it is easier to be implemented
than  CNFD, ReFD and SIFD.
It is an implicit scheme  where, at each time step, only a linear system (whose coefficient matrix
is time independent) needs to be solved,
which can be done in $\mathcal{O}(J)$ operations in 1D and
 $\mathcal{O}(J^d\ln J)$ operations in 2D/3D. Therefore it is significantly cheaper than  the
CNFD and ReFD methods. The TSFD method is second-order accurate both in
time and space for any fixed $\varepsilon=\varepsilon_0=\mathcal{O}(1)$
\cite{BaoCai313,BaoTangXu}.  In addition, it conserves the {\sl mass}
(\ref{massConserveC1d}) at the discretized level as (\ref{massdis})
\cite{BaoCai313,BaoTang,BaoTangXu} and thus it is unconditionally stable;  it is also {\sl time transverse invariant}.
However, there is no  {\sl energy conservation}
property of the NLSE/GPE (\ref{NLSE1}) at the discretized level
\cite{BaoCai313,BaoTang,BaoTangXu}.

For comparison and convenience of the readers,  we summarize the main physical and computational
properties of the above popular numerical methods CNFD, ReFD, SIFD, TSSP and TSFD
in Table \ref{SchemeProperties}. From this Table, we can see that the TSSP method shares the most
properties as the original NLSE/GPE. In addition, it is very efficient and accurate as well as
easy to be implemented in practical computations, especially in 2D and 3D \cite{BaoCai313,BaoJakschP,BaoJinP,BaoJinP2,BaoShen,Thalhammer}. In fact, the TSSP method
becomes more and more popular in practical computations, especially in the numerical simulation of
Bose-Einstein condensation \cite{BaoCai313,BaoJakschP,BaoShen,Huang,Huang1,Thalhammer,Thalhammer2}.

\begin{savenotes}
  \begin{table}[htbp]
    \centering
    \begin{tabular}{|c||c|c|c|c|c|c|r||}
      \hline
      Method  & TSSP  & CNFD  & SIFD  &  ReFD &TSFD\\
      \hline
      Time Reversible & Yes & Yes &  Yes &  Yes & Yes\\
      Time Transverse Invariant & Yes & No &  No  & No & Yes\\
      Mass Conservation & Yes & Yes &  No &  Yes & Yes\\
      Energy Conservation & No & Yes &  No &  Yes\footnote{Only for cubic nonlinearity} & No\\
      Dispersion Relation & Yes & No &  No &  No  & Yes\\

      \hline
      Unconditional Stability & Yes & Yes &  No  & Yes  & Yes\\
      Explicit  Scheme & Yes & No &  No &  No & No \\
      Time Accuracy & $2^{\textrm{th}}$ or $4^{\textrm{th}}$ & $2^{\textrm{th}}$ &  $2^{\textrm{th}}$ &  $2^{\textrm{th}}$  & $2^{\textrm{th}}$\\
      Spatial Accuracy & spectral & $2^{\textrm{th}}$ & $2^{\textrm{th}}$  &  $2^{\textrm{th}}$ &  $2^{\textrm{th}}$\\
      Memory Cost & $\mathcal{O}(J^d)$ & $\mathcal{O}(J^d)$ &  $\mathcal{O}(J^d)$ &  $\mathcal{O}(J^d)$ & $\mathcal{O}(J^d)$ \\
      Computational Cost & $\mathcal{O}(J^d\log J)$ & $\gg\mathcal{O}(J^d)$\,\footnote{Depends on the solver for the nonlinear system} &  $\mathcal{O}(J^d\log J)$\,\footnote{If $d=1$, $\mathcal{O}(J)$} &  $\mathcal{O}(J^d\log J)$\,\footnote{If $d=1$, $\mathcal{O}(J)$} & $\mathcal{O}(J^d\log J)$\,\footnote{If $d=1$, $\mathcal{O}(J)$} \\ \hline
      Resolution when $0<\varepsilon \ll 1$\,\footnote{For cubic repulsive nonlinearity} &
      $\begin{array}{c} h=\mathcal{O}(\varepsilon)\\
      \tau=\mathcal{O}(\varepsilon)
      \end{array}$ & $\begin{array}{c} h=o(\varepsilon)\\
      \tau=o(\varepsilon) \end{array}$ & $\begin{array}{c} h=o(\varepsilon)\\
      \tau=o(\varepsilon) \end{array}$ & $\begin{array}{c} h=o(\varepsilon)\\
      \tau=o(\varepsilon)
      \end{array}$ & $\begin{array}{c} h=o(\varepsilon)\\
      \tau=o(\varepsilon)
      \end{array}$\\
      \hline
    \end{tabular}
    \caption{Physical and numerical properties of different popular numerical methods in the $d$-dimensional case
    with $J$  unknowns in each direction.\label{SchemeProperties}}
  \end{table}
\end{savenotes}

\noindent {\sl Remark 2.1}
If the homogeneous Dirichlet boundary condition (\ref{NLSE2}) is replaced by the periodic boundary condition, e.g.
for BEC on a ring \cite{BaoCai313,BaoJinP,BaoJinP2,BaoTangXu},
the above numerical methods are still valid provided that we
replace $1\le j\le J-1$ by $0\le j\le J-1$ in all the numerical methods, and replace
$\psi_0^{n+1}=\psi_J^{n+1}=0$ by $\psi_0^{n+1}=\psi_J^{n+1}$ and $\psi_{-1}^{n+1}=\psi_{J-1}^{n+1}$
in the CNFD, ReFD and SIFD methods,  the sine basis by the Fourier basis in the TSSP method,
and  $\psi_0^{(2)}=\psi_J^{(2)}=0$ by $\psi_0^{(2)}=\psi_J^{(2)}$ and $\psi_{-1}^{(2)}=\psi_{J-1}^{(2)}$
in the TSFD method, respectively. Similarly, if the homogeneous Dirichlet boundary condition (\ref{NLSE2})
is replaced by the homogeneous Neumann boundary condition, e.g.
for the dynamics of dark solitons and their interaction \cite{Bao2,BaoCai313,BaoTangXu},
the above numerical methods are still
valid provided that we replace $1\le j\le J-1$ by $0\le j\le J$ in all the numerical methods, and replace
$\psi_0^{n+1}=\psi_J^{n+1}=0$ by $\psi_{-1}^{n+1}=\psi_0^{n+1}$ and $\psi_{J+1}^{n+1}=\psi_{J-1}^{n+1}$
in the CNFD, ReFD and SIFD methods,  the sine basis by the cosine basis in the TSSP method,
and  $\psi_0^{(2)}=\psi_J^{(2)}=0$ by $\psi_{-1}^{(2)}=\psi_0^{(2)}$ and $\psi_{J+1}^{(2)}=\psi_{J-1}^{(2)}$
in the TSFD method, respectively
 \cite{Bao2,BaoCai313,BaoTangXu}. Then the physical and numerical properties of these numerical
methods are still valid.

\subsection{Other numerical methods}

In the literatures, there are many other numerical methods proposed for discretizing the NLSE/GPE.
For example, Sanz-Serna  \cite{sanzserna,RFH93,Sanz1984,SV86} proposed the following second-order
finite difference (SSFD) method which is well-adapted to the computation of soliton-like solutions
of NLSE/GPE
\begin{equation}
i\varepsilon \frac{\psi^{n+1}_j-\psi^{n}_j}{\tau} =- \frac{\varepsilon^2}{4}\left[\left(D^2_x\boldsymbol{\psi}^{n+1}\right)_j+
\left(D^2_x\boldsymbol{\psi}^n\right)_j\right]
+ \left[V(x_j)
+f\left(\left|\frac{\psi^{n+1}_j+\psi^{n}_j}{2}\right|^{2}\right)\right] \frac{\psi^{n+1}_j+\psi^{n}_j}{2},
\quad 1\le j\le J-1, \ n\ge0,
\label{SSsemidiscrete}
\end{equation}
and the boundary and initial conditions (\ref{NLSE2})-(\ref{NLSE3}) are
discretized as (\ref{FDmatrix}). Another one is the leap-frog finite difference (LPFD) method as
\cite{MurA,ZhangBao}
\begin{equation}
i\varepsilon \frac{\psi^{n+1}_j-\psi^{n-1}_j}{2\tau} =- \frac{\varepsilon^2}{2}\left(D^2_x\boldsymbol{\psi}^{n}\right)_j
+ \left[V(x_j)
+f\left(\left|\psi^{n}_j\right|^{2}\right)\right] \psi^{n}_j,
\qquad 1\le j\le J-1, \qquad n\ge1,
\label{lpidiscrete}
\end{equation}
and the boundary and initial conditions (\ref{NLSE2})-(\ref{NLSE3}) are
discretized as (\ref{FDmatrix}) and the first step can be computed as (\ref{fstep87}).
 In fact, in the physics literatures,
for time discretization,  the fourth-order Runge-Kutta (RK4) method was also used \cite{Adh1,Baer,Caradoc,CD2},
which is not time symmetric. Thus, in general, it should be avoid to use RK4 for solving
the NLSE/GPE. In some mathematics literatures,
for space discretization, the 4th-order compact finite difference method has been used \cite{Liao,Xie}
and  the finite element (FE) method was also developed and analyzed
\cite{AD91,ADK91}. However, since the computational domains for NLSE/GPE are usually simple and regular,
the FE method for spatial discretization is usually not adapted in practical computations.
Last but not the least, we want to mention that, in some physics papers,
for solving  NLSE/GPE in 2D or 3D, the alternating direction implicit (ADI)
is first adapted to decouple the NLSE/GPE into dimension-by-dimension and then the SIFD or
CNFD is used to discretize each NLSE/GPE in 1D \cite{Adh1,Edwards,Gao,MurA,Ruprecht,Saito}. For other more
numerical methods for NLSE/GPE, we refer to
\cite{A93,Arnold1,Baer,Bao2,Tosi2,Cerim,Chan2,Chan1,Chang2,Delbour,Dion,Dorfler,Fei,Tomio,
Griff,Hong,Hong3,Hong1,Ignat,Ismail1,Ismail2,Kyza,MurA,RFH93,Ruprecht,Wu} and references therein.
In general, these numerical methods have less good properties
as those methods mentioned in the previous subsection and thus we omit
the details here for brevity.

\subsection{Perfectly matched layers and/or absorbing boundary conditions}\label{ABC}

In some situations for solving the NLSE/GPE (\ref{eqSchrodinger}) numerically,
e.g. very long time dynamics and the potential is not a
confinement potential in nonlinear optics \cite{Arnold01,Arnold},
the simple boundary condition, such as homogeneous Dirichlet or Neumann or periodic
boundary condition used in the previous subsections to truncate (or approximate)
the original NLSE/GPE from the whole space problem
to a bounded computational domain,  might bring large truncation errors
except that the bounded computational domain
is chosen extremely large and/or time-dependent. Thus, in order to choose a smaller computational domain
which might save memory and/or computational cost, perfectly matched layers (PMLs) \cite{BerengerMPL} or
high-order absorbing (or artificial) boundary conditions (ABCs)
\cite{Review08,Arnold01,Arnold,BaylissTurkel,EngquistMajda,Mur,Pang}
need to be designed and/or used at the artificial boundary so that one can
truncate (or approximate) the original NLSE/GPE into a smaller bounded computational domain.
Over the last 20 years, different PMLs \cite{NissenKreiss,ZhengPML}
and/or ABCs \cite{AntBesDes06,ABK1dP,ABPLaser,ABKNLS2D,ABK2DPTheory,ABK2DPNumerics,Arnold01,Arnold,
SzeftelNLS,Sze04,ZhengCubic}
have been designed for solving the NLSE/GPE in the literatures.
Here we simply review some of them for the completeness of this paper.

In fact, PMLs were introduced by B\'erenger in 1994 for electromagnetic field computations \cite{BerengerMPL}
and they have been extended to NLSE/GPE by various authors recently \cite{NissenKreiss,ZhengPML}.
Again, here we present the idea in 1D. Suppose that  one is only interested in the solution of the
NLSE/GPE (\ref{eqSchrodinger}) with $d=1$ over the physical domain $(a,b)$. One introduces two layers with width
$R_0>0$ at $x=a$ and $x=b$ and defines  the following function
\[S(x):=1+R\sigma(x), \qquad \tilde{a}:=a-R_0\le x\le b+R_0:=\tilde{b},\]
where $R\in {\mathbb C}$ and $\sigma(x)$ is
a real-valued function which must be chosen properly and very carefully.
The goal is to artificially damp the wave traveling inside the two layers to zero
without modifying their dynamics in $(a,b)$ by properly chosen $R$ and $\sigma(x)$ \cite{Review08,Arnold01,NissenKreiss,ZhengPML}.
Different choices have been proposed in the literatures
\cite{Review08,NissenKreiss,ZhengPML}. Among them, the following choice
works well for the NLSE/GPE (\ref{eqSchrodinger}) in 1D with cubic focusing nonlinearity
\cite{Review08,NissenKreiss,ZhengPML}
 \begin{equation}
 R=e^{i\pi/4},\qquad  \qquad \qquad \label{eq:sigma-PML-V}
  \sigma(x) =
  \frac{1}{\delta^{2}}\begin{cases}
    (x-a)^2, &\qquad \tilde{a}\le x<a , \\
    0, &\qquad a\le x \le b, \\
    (x-b)^2 , &\qquad b<x \le \tilde{b},
  \end{cases}
\end{equation}
where  $\delta$ is a positive constant. Then, the NLSE/GPE (\ref{eqSchrodinger})
with $d=1$ will be truncated (or approximated) as
\begin{eqnarray}\label{NLSE176}
&&i \varepsilon \partial_t \psi(t,x) = -\frac{\varepsilon^2}{2} \frac{1}{S(x)}\partial_x\left(\frac{1}{S(x)}\partial_x\psi(t,x)\right)
+V(x) \psi(t,x)  +f(|\psi(t,x) |^2) \psi(t,x) , \qquad \tilde{a}<x<\tilde{b},\quad t>0,\qquad \\
\label{NLSE276}
&&\psi(t,\tilde{a})=\psi(t,\tilde{b})=0, \qquad t\ge0,\\
 \label{NLSE376}
&&\psi(t=0,x)=\psi_0(x),\qquad \tilde{a}\le x\le \tilde{b}.
\end{eqnarray}
The above problem can be discretized by CNFD, SIFD, ReFD and TSFD methods straightforward
provided that we introduce the following finite difference operator to approximate
$\left. \frac{1}{S(x)}\partial_x\left(\frac{1}{S(x)}\partial_x d(x)\right)\right|_{x=x_j}$ as
 $$
\left. \frac{1}{S(x)}\partial_x\left(\frac{1}{S(x)}\partial_x d(x)\right)\right|_{x=x_j}\approx
 \frac{1}{2h^{2}S_{j}}\left[\frac{1}{S_{j-1/2}}d_{j-1} -
 \left(\frac{1}{ S_{j-1/2}}+\frac{1}{S_{j+1/2}}\right)d_{j} +
\frac{1}{S_{j+1/2}}d_{j+1}\right],
$$
where $S_j=S(x_j)$, $S_{j-1/2}=S(x_j-h/2)$ and $d_j$ is the approximation of $d(x_j)$.
This PML can be  efficient and accurate for the NLSE/GPE in 1D in some cases if one can choose
 $\delta$ and $R_0$ properly. Of course, extensions to 2D and 3D are a little bit more difficult
and it might give bad numerical results or even  not work in some situations when one
chooses improper layer function $S(x)$ and/or layer width $R_0$.
Another way to design the PML is to introduce a damping potential which is also called
as complex absorbing potential (CAP) or exterior complex scaling (ECS) \cite{Juengel,Scrinzi},
i.e. choosing $S(x)\equiv1$ for $\tilde{a}\le x\le \tilde{b}$ and
replacing the real-valued potential $V(x)$ by a
 complex-valued  potential $\tilde{V}(x):=V(x)-i\sigma(x)$,
 where $\sigma(x)$ is given in (\ref{eq:sigma-PML-V})
with $\delta$ and $\sigma_0$ two real constants to be determined based on the application.
Then, the problem can be discretized by TSSP, CNFD, SIFD, ReFD and TSFD methods straightforwardly.
Again, this CAP or ECS can be  efficient and accurate for the NLSE/GPE in 1D in some situations if one can choose
 $\delta$ and $R_0$ suitably. The extensions to  the 2D and 3D cases are direct.
 It might also lead to large inaccuracies in the calculations, more particularly for nonlinear problems.
 In general, PMLs are more robust from the accuracy point of view \cite{ABPLaser}.

Another way is to find the Dirichlet-to-Neumann  (DtN) map for the Schr\"{o}dinger operator
without and/or with external potential and nonlinearity at the boundary of the bounded
computational domain by using the continuous and/or discrete Laplace transform
\cite{Review08,AntBes03,AntBesDes06,ABK1dP,ABKNLS2D,ABK2DPTheory,ABK2DPNumerics,Arnold01,Parabolic}.
The DtN map is usually nonlocal and a time-dependent pseudo-differential operator, and thus
its approximations which involve time fractional derivatives and integrals,
are usually used in practical computations. These boundary conditions are also called as
absorbing (or artificial or transparent) boundary conditions according to their properties.
The goal is to avoid or at least to minimize the wave
reflected back inside the domain while it should be outgoing.
Since there are several very nice review papers on this part,
we omit the details here for brevity and refer to  \cite{Review08}
and references therein. Due to the nonlocal ABCs at the artificial boundary,
one can usually solve the problem on the bounded computational domain
by the second-order finite difference, such as ReFD or SIFD method.
It is usually very hard to design spectral order method in space
for the NLSE/GPE with nonlocal ABCs.

For comparison of the performance and effectiveness of
different PMLs and/or ABCs for NLSE/GPE and their applications, we refer to
\cite{ABKNLS2D,PABE,XuHan2006,XuHanWu2007,ZhangZhenliWu2009} and references therein.

\subsection{Numerical methods in the semi-classical regime}
When $0<\varepsilon \ll1$ in the NLSE/GPE (\ref{eqSchrodinger}), i.e.
in the semiclassical limit regime, then the solution propagates waves
with wave-length at $\mathcal{O}(\varepsilon)$ in both space and time
\cite{BaoJinP,BaoJinP2,Carles01,GMMP,Huang,Jin}.
The highly oscillatory nature brings severe numerical burdens
in numerical computation of the solution of the NLSE/GPE (\ref{eqSchrodinger}) \cite{Cheng,Huang,Jin}.
In fact, in order to capture numerically ''correct'' physical observables such as density and current,
different numerical methods request different mesh strategies (or $\varepsilon$-scalability)
for discretizing the NLSE/GPE (\ref{eqSchrodinger}) in the semiclassical limit regime,
i.e. $0<\varepsilon \ll1$. Based on the analysis and extensive numerical studies
\cite{BaoJinP,BaoJinP2,Huang,Jin},
it has been found that the $\varepsilon$-scalability for TSSP is
mesh size $h=\mathcal{O}(\varepsilon)$, and time step
$\tau=\mathcal{O}(1)$-independent of $\varepsilon$ and
 $\tau=\mathcal{O}(\varepsilon)$ for the linear Schr\"{o}dinger equation and NLSE/GPE,
respectively;  for the CNFD, ReFD, SIFD and TSFD as well as many other finite difference methods,
the mesh size is $h=o(\varepsilon)$ and the time step is
 $\tau=o(\varepsilon)$ \cite{BaoJinP,BaoJinP2,Jin,MaPiPo,MPPS}. Thus, among those numerical methods for discretizing
 the NLSE/GPE (\ref{eqSchrodinger}) directly, the TSSP method has the best resolution
 in the semiclassical limit regime.

Of course, there are many other efficient and accurate numerical methods for
solving the NLSE/GPE (\ref{eqSchrodinger}) based on the oscillatory structure
of the solution. For the linear Sch\"{o}dinger equation, the Gaussian beam
method \cite{Faou0,Jin,Jin23,Jin24,Leung,Qian}
shows good resolution in spatial discretization, which in general
requires  $h=\mathcal{O}(\sqrt{\varepsilon})$ and
$\tau=\mathcal{O}(1)$-independent of $\varepsilon$.
For the details of the Gaussian beam method, we refer to
a recent nice review paper \cite{Jin} and references therein.
However, in general, the Gaussian beam method cannot be
extended to deal with the NLSE/GPE (\ref{eqSchrodinger})
due to the fact that the superposition is used in the method. We also
mention here that there are some numerical methods in the literatures
for solving the linear Schr\"{o}dinger equation based on the Liouville equation through
the  Wigner transform, which can give good results in the linear case in the semiclassical limit
regime \cite{Cheng,GMMP,Jin}.

Another approach for dealing with the nonlinear (or linear)
Schr\"{o}dinger equation (\ref{eqSchrodinger}) in the semiclassical limit regime
is through the Madelung (or WKB) expansion \cite{Mad,Carles01}.
In the semiclassical limit regime, i.e. $0<\varepsilon\ll1$,
we denote the solution $\psi$ to the NLSE/GPE (\ref{eqSchrodinger}) as $\psi^\varepsilon$
and use the following WKB ansatz (or Madelung transformation) \cite{Carles01,Jin,Mad}
\begin{equation}\label{wkb1}
\psi^\varepsilon(t,\mathbf{x})=\sqrt{\rho^\varepsilon(t,\mathbf{x})}\;\exp\left({\frac{i\,
  S^\varepsilon(t,\mathbf(x)}{\varepsilon}}\right), \qquad \xbf\in {\mathbb R}^d,\quad t\ge0,
\end{equation}
where the real-valued functions $\rho^\varepsilon=|\psi^\varepsilon|^2$ and $S^\varepsilon$
are the density and phase, respectively.
Plugging the above WKB ansatz into the NLSE/GPE (\ref{eqSchrodinger}) and identifying real and imaginary
parts, the NLSE/GPE is reformulated as a coupled system for density
 and quantum velocity
$v^\varepsilon=\nabla S^\varepsilon$, which is known as quantum hydrodynamic (QHD) system
(made of a compressible, isentropic Euler system with Bohm
potential) \cite{Carles01,Carles02}. This model was used in \cite{Arnold1,Carles02,Cheng,DGM,Jin22} to build an
asymptotic preserving (AP) scheme which allows to compute numerical
solution with time step and mesh size independent of
$\varepsilon$. However, the QHD system fails to represent
the solution of the original NLSE/GPE near vacuum, i.e. $\rho^\varepsilon=0$
\cite{Carles01,Carles02}, and thus the scheme in \cite{DGM} suffers from difficulties when the
density $\rho^\varepsilon$ vanishes in the domain
and/or shocks or sharp changes happen in the QHD system.
To overcome this drawback in the Madelung (or WKB) expansion for including vacuum,
Grenier \cite{Grenier98} introduced a modified Madelung (or WKB) expansion
with the following ansatz  \cite{Grenier98,Carles01}
\begin{equation}\label{wkb123}
\psi^\varepsilon(t,\mathbf{x})=A^\varepsilon(t,\xbf)\;\exp\left({\frac{i\,
  S^\varepsilon(t,\mathbf(x)}{\varepsilon}}\right), \qquad \xbf\in {\mathbb R}^d,\quad t\ge0,
\end{equation}
where $A^\varepsilon:=a^\varepsilon+i\,b^\varepsilon$ is a
complex-valued function with $a^\varepsilon$ and $b^\varepsilon$ two real-valued functions.
Plugging (\ref{wkb123})
into the NLSE/GPE (\ref{eqSchrodinger}) and collecting $O(1)$ and $O(\varepsilon)$-terms,
one obtains a system for $A^\varepsilon$ and $S^\varepsilon$ (or $v^\varepsilon=\nabla
S^\varepsilon$) \cite{Carles01,Grenier98}. Based on this formulation, an AP scheme
was proposed in \cite{remibijan}, which works up to the time when caustics
happens. More recently, another AP scheme \cite{BCM2013,Jin22}
was presented based on the Grenier's expansion with a regularized term,
which can work after the caustics happens. For more details, we
refer to \cite{BCM2013,remibijan} and references therein.

\section{Extension to NLSE/GPE with damping and angular rotation terms} \label{s3}
\setcounter{equation}{0}

\subsection{For damped NLSE/GPE}

As proven in \cite{BaoCai313,Cazenave,SulemSulem}, finite time blow-up may happen
for the NLSE/GPE (\ref{eqSchrodinger}) with a focusing nonlinearity,
e.g. the cubic nonlinearity (\ref{nls_cubic}) with $\beta<0$ in 2D/3D.
However, the physical quantities modeled by $\psi:=\psi(t,\xbf)$ do not become
infinite, e.g. in BEC, which implies that the validity of (\ref{eqSchrodinger}) breaks
down near the singularity. Additional physical mechanisms, which
were initially small, become important near the singular point and
prevent the formation of the singularity \cite{Donley,Fibich,FibichP}. In BEC, the particle
density $\rho:=|\psi|^2$ becomes large close to the critical point and
inelastic collisions between particles which are negligible for
small densities become important \cite{BaoJaksch,Donley,Fibich,FibichP}. Therefore, a small damping
(absorption) term is introduced into the NLSE/GPE (\ref{eqSchrodinger}) which
describes inelastic processes \cite{BaoCai313,BaoJaksch,BaoJakschP46,Donley,Fibich,FibichP}. We are interested in the cases
where these damping mechanisms are important and, therefore,
restrict ourselves to the case of focusing nonlinearity, i.e. $\beta<0$,
where $\beta$ may also be time dependent. We consider the following
damped NLSE/GPE \cite{BaoCai313,BaoJaksch,BaoJakschP46,Fibich,FibichP}
\begin{eqnarray} \label{eq:sdged:sec2}
&&i\; \partial_t\psi(t,\xbf)=-\frac{1}{2}\;\nabla^2 \psi(t,\xbf)+ V(\xbf)\; \psi(t,\xbf)
+f(|\psi(t,\xbf)|^2)\psi(t,\xbf)-i\; g(|\psi(t,\xbf)|^2)\psi(t,\xbf),
\qquad t>0, \quad \xbf\in {\mathbb R}^d, \qquad \quad\\
\label{eq:sdgid:sec2}
&&\psi(t=0,\xbf)=\psi_0(\xbf), \qquad \xbf\in {\mathbb R}^d,
\end{eqnarray}
where $g(\rho)\ge 0$ for $\rho:=|\psi|^2\ge 0$ is a real-valued
monotonically increasing function.

The general form of (\ref{eq:sdged:sec2}) covers many damped NLSE/GPE arising
in various different applications. In BEC, for example, when
$g(\rho)\equiv 0$, (\ref{eq:sdged:sec2}) reduces to the usual NLSE/GPE
(\ref{eqSchrodinger}); a linear damping term $g(\rho)\equiv \delta$ with
$\delta>0$ describes inelastic collisions with the background gas;
cubic damping $g(\rho)=\delta_1 \rho$ with $\delta_1>0$ corresponds
to two-body loss \cite{Saito}; and a quintic damping term
of the form $g(\rho)=\delta_2  \rho^2$ with $\delta_2>0$ adds
three-body loss to the NLSE/GPE (\ref{eqSchrodinger}) \cite{Saito}. It is
easy to see that the decay of the {\sl mass} according to
(\ref{eq:sdged:sec2}) due to damping is given by
\begin{equation} \label{eq:normNt:sec2}
\dot{N}(t)=\frac{d}{ dt} \int_{{\mathbb R}^d}\; |\psi(t,{\bf x})|^2\;
d{\bf x} = -2\int_{{\mathbb R}^d}\; g(|\psi(t,{\bf x})|^2)|\psi(t,{\bf
x})|^2\; d{\bf x} \le 0, \qquad t >0.
\end{equation}
Particularly, if
$g(\rho)\equiv \delta$ with $\delta>0$, the mass is given by
\begin{equation} \label{eq:drnt:sec2}
N(t)=\int_{{\mathbb R}^d}\; |\psi(t,{\bf x})|^2\, d{\bf x}= e^{-2\delta\, t} N(0)=
e^{-2\delta\, t}\; \int_{{\mathbb R}^d}\; |\psi_0({\bf x})|^2\, d{\bf x}, \qquad t\ge 0.
\end{equation}

Due to the appearance of the damping term, new ideas are needed to deal with them
and different numerical methods have been presented in the literatures \cite{BaoJaksch,BaoJakschP46}.
In fact, the numerical methods such as TSSP, CNFD, ReFD, SIFD and TSFD methods for the
NLSE/GPE (\ref{eqSchrodinger}) presented in the previous section
can be easily extended to the damped NLSE/GPE (\ref{eq:sdged:sec2}).
For simplicity of notations, here we only present the TSSP method
for (\ref{eq:sdged:sec2}) with quintic damping term in 1D, i.e.
$d=1$ and $g(\rho)=\delta_2  \rho^2$ with $\delta_2>0$.
From time $t=t_n$ to time $t=t_{n+1}$, the damped NLSE/GPE (\ref{eq:sdged:sec2}) is solved in two steps.
One solves
\begin{equation} \label{eq:1step:sec4}
i\, \partial_t\psi(t,x)=- \frac{1}{2} \partial_{xx} \psi(t,x), \qquad a<x<b, \quad t>t_n,
\end{equation}
with homogeneous Dirichlet boundary condition for one time step of length $\tau$, followed by solving
\begin{equation}
\label{eq:2step:sec4} i\, \partial_t\psi(t,x)= V(x)\psi(t,x)+
f(|\psi(t,x)|^{2}) \psi(t,x)-i\,\delta_2|\psi(t,x)|^4\psi(t,x), \qquad a\le x\le b, \qquad t>t_n,
\end{equation}
 for the same time step $\tau$. Again, Eq. (\ref{eq:1step:sec4}) is
discretized in space by the sine-spectral method and integrated in
time {\it exactly}. For $t\in[t_n,t_{n+1}]$, multiplying the ODE
(\ref{eq:2step:sec4}) by $\overline{\psi(t,x)}$ and then
subtracting from its  conjugate,  we obtain for $\rho(t,x):=|\psi(t,x)|^2$ \cite{BaoJaksch}
\begin{equation}
\label{eq:sstepc:sec4}
\partial_t \,\rho(t,x)=-2 \delta_2\, \rho^3(t,x), \qquad t> t_n, \qquad a\le x\le b,
\end{equation}
which can be solved analytically as
\begin{equation}\label{rhoexact}
\rho(t,x)=\frac{\rho(t_n,x)}{\sqrt{4 \delta_2 (t-t_n)+\rho^2(t_n,x)}}, \qquad t\ge t_n,
\quad a\le x\le b.
\end{equation}
Plugging (\ref{rhoexact}) into (\ref{eq:2step:sec4}),
we get a linear ODE as
\begin{equation}
\label{eq:2step:sec496} i\, \partial_t\psi(t,x)=\left[ V(x)+
f\left(\frac{\rho(t_n,x)}{\sqrt{4 \delta_2 (t-t_n)+\rho^2(t_n,x)}}\right)
-i\,\delta_2\frac{\rho^2(t_n,x)}{4 \delta_2 (t-t_n)+\rho^2(t_n,x)}\right]\psi(t,x),
\quad t>t_n, \ \ a\le x\le b,
\end{equation}
which can be integrated {\sl exactly} as for $0\le s\le \tau$ and $a\le x\le b$
\begin{equation}
\psi(t_n+s,x)=\psi(t_n,x)\; \exp\left[-i\left(V(x)s+\frac{i}{4}\rho^2(t_n,x)\ln\frac{\rho^2(t_n,x)}
{4\delta_2s+\rho^2(t_n,x)}+\int_0^{s}f\left(\frac{\rho(t_n,x)}{\sqrt{4 \delta_2 u+\rho^2(t_n,x)}}\right)du\right)\right].
\end{equation}
For cubic nonlinearity, the last term in the above equation can be integrated {\sl analytically}.
For a more general nonlinearity, if it cannot be integrated analytically, one can use
a numerical quadrature, e.g. the Simpson's rule, to evaluate it numerically \cite{BaoJaksch,BaoJakschP46}.
Then, we can construct the second-order time-splitting sine pseudospectral (TSSP)
 method for the damped NLSE/GPE (\ref{eq:sdged:sec2}) \textit{via} the Strang
 splitting \cite{BaoJaksch,BaoJakschP46}; the details are omitted here for brevity.

\subsection{NLSE/GPE with an angular rotation term}

In view of potential applications of BEC,
the study of quantized vortices, which are
related to superfluid properties,
is one of the key issues.
Currently, one of the most popular ways to generate quantized vortices from
BEC ground state is the following:  impose a laser beam rotating with an angular
velocity
on the magnetic trap holding the atoms to
create a harmonic anisotropic potential. Various experiments have confirmed the observation of quantized vortices in  BEC under a rotational frame \cite{Abo,Aft,AftalionDu,BaoCai313,Caradoc,Madison}.
At temperatures $T$ much smaller than the critical temperature
$T_c$,  BEC in a rotational frame  is well described by
the macroscopic wave function $\psi:=\psi(t,\xbf)$, whose evolution is
governed by the following dimensionless GPE
with an angular momentum  term around the $z$-axis
\cite{Aft,AftalionDu,BaoCai313,BaoShen,BaoWang,Fetter,PitaevskiiStringari}:
\begin{equation} \label{eq:gperot:sec5}
i\,\partial_t \psi(t,\xbf)=-\frac{1}{2}\nabla^2 \psi(t,\xbf) +
V(\xbf) \psi(t,\xbf)+ \beta |\psi(t,\xbf)|^2\psi(t,\xbf)-\Omega
L_z\,\psi(t,\xbf), \qquad \xbf\in{\mathbb R}^d, \quad t>0,
\end{equation}
with initial data
\begin{equation} \label{initgpe876}
\psi(t=0,\xbf)=\psi_0(\xbf), \qquad \xbf\in{\mathbb R}^d,
\end{equation}
where $d=2$ or $3$ for 2D and 3D, respectively,
$\beta$ is a dimensionless constant describing the interaction strength,
$V:=V(\xbf)$ is a given real-valued potential which is usually chosen as a
harmonic potential,
$\Omega$ is the dimensionless rotation velocity, and
\begin{equation}\label{eq:rota:sec5}
L_z=-i \left(x\partial_y -y\partial_x\right)
\end{equation}
is the $z$-component of the angular momentum operator $\boldsymbol{L}=(L_x,L_y,L_z)^T$ given by
$\boldsymbol{L}=\xbf\wedge \boldsymbol{P}$, with the momentum $\boldsymbol{P}=-i\,\nabla$.
The appearance of the angular momentum term
means that we are using a reference frame where the trap is at rest.
The above GPE is {\sl time reversible} and {\sl time transverse invariant} and
it conserves the {\sl mass} (\ref{massConserveC}) and the {\sl energy} defined as
\cite{Aft,AftalionDu,BaoCai313,BaoShen,BaoWang,Fetter,PitaevskiiStringari}
\begin{eqnarray}
E(t):=\int_{\mathbb{R}^{d}} \left[\frac{1}{2} |\nabla \psi(t,\xbf) |^{2} + V(\xbf) | \psi(t,\xbf)|^{2} +\frac{\beta}{2}|\psi(t,\xbf)|^4-\Omega \overline{\psi(t,\xbf)}\;L_z\psi(t,\xbf)
\right] d \mathbf{x} \equiv E(0), \qquad  t \ge 0.
\label{energyConserveC123}
\end{eqnarray}

Due to the appearance of the angular momentum term, new difficulties are introduced
in solving the GPE (\ref{eq:gperot:sec5}) for rotating BEC numerically.
Several efficient and accurate numerical methods have been proposed in the
 literatures for discretizing it \cite{BaoCai313,BaoCai2,BaoDuZhang,BaoLiShen,BaoMTZ,BaoWang,Caradoc,Saito}.
In fact, the numerical methods such as TSSP, CNFD, ReFD, SIFD and TSFD methods for the
NLSE/GPE (\ref{eqSchrodinger}) presented in the previous section
can be easily extended to the GPE (\ref{eq:gperot:sec5}) with an angular momentum rotation.
For conciseness,  we only present here the TSSP method
for (\ref{eq:gperot:sec5}) in 2D, i.e. $d=2$.

From time $t=t_n$ to $t=t_{n+1}$, the GPE (\ref{eq:gperot:sec5}) is solved in two steps.
One solves
\begin{equation} \label{eq:1step:sec4756}
i\, \partial_t\psi(t,\xbf)=- \frac{1}{2} \nabla \psi(t,\xbf)-\Omega
L_z\,\psi(t,\xbf), \qquad  t>t_n,
\end{equation}
for one time step of length $\tau$, followed by solving the ODE
\begin{equation}
\label{eq:2step:sec497} i\, \partial_t\psi(t,\xbf)= V(x)\psi(t,\xbf)+
\beta|\psi(t,\xbf)|^2 \psi(t,\xbf),  \qquad t>t_n,
\end{equation}
for the same time step $\tau$. Similar to (\ref{NLSE2nd}), Eq. (\ref{eq:2step:sec497})
can be solved {\sl analytically} \cite{BaoCai313,BaoDuZhang}. Some numerical methods
have been presented for discretizing (\ref{eq:1step:sec4756}).
One numerical method is to adapt the polar coordinates $(r,\theta)$
in 2D such that the angular momentum
rotation becomes constant coefficient, and then to discretize it in  the
$\theta$-direction by the Fourier spectral method, in the $r$-direction by the second-order
or fourth-order finite difference method or finite element method,
and in time by the Crank-Nicolson method. For more details, we refer to
\cite{BaoCai313,BaoDuZhang}. Another method is to apply the Alternating Direction Implicit (ADI) method
to decouple  (\ref{eq:1step:sec4756}) into two sub-problems as
\begin{eqnarray} \label{eq:1step:sec411}
&&i\, \partial_t\psi(t,\xbf)=- \frac{1}{2} \partial_{xx} \psi(t,\xbf)-i\Omega y\partial_x
\,\psi(t,\xbf), \qquad  t>t_n, \\
\label{eq:1step:sec412}
&&i\, \partial_t\psi(t,\xbf)=- \frac{1}{2} \partial_{yy} \psi(t,\xbf)+i\Omega x\partial_y
\,\psi(t,\xbf), \qquad  t>t_n.
\end{eqnarray}
Now the first problem  (\ref{eq:1step:sec411}) is constant coefficient with respect to $x$
and the second problem (\ref{eq:1step:sec411}) is constant coefficient with respect to $y$,
and thus they can be discretized in space by the Fourier spectral method and
integrated in time {\sl exactly}. Again, for more details, we refer to
\cite{BaoCai313,BaoWang}.

A different way to apply the time-splitting technique to the
GPE (\ref{eq:gperot:sec5}) is the following: from time $t=t_n$ to time $t=t_{n+1}$,
one solves
\begin{equation} \label{eq:1step:sec445}
i\, \partial_t\psi(t,\xbf)=- \frac{1}{2} \nabla^2 \psi(t,\xbf)+\frac{|\xbf|^2}{2}\psi(t,\xbf)-\Omega
L_z\,\psi(t,\xbf), \qquad \xbf\in{\mathbb R}^2,\quad  t>t_n,
\end{equation}
for one time step of length $\tau$, followed by solving
\begin{equation}
\label{eq:2step:sec446} i\, \partial_t\psi(t,\xbf)= \left(V(x)-\frac{|\xbf|^2}{2}\right)\psi(t,\xbf)+
\beta|\psi(t,\xbf)|^2 \psi(t,\xbf),  \qquad \xbf\in{\mathbb R}^2,\quad t>t_n,
\end{equation}
for the same time step $\tau$. Again, similar to (\ref{NLSE2nd}), Eq. (\ref{eq:2step:sec446})
can be solve {\sl analytically} \cite{BaoCai313,BaoLiShen}. Eq. (\ref{eq:1step:sec445})
can be discretized in space by the generalized-Laguerre-Fourier spectral method and
integrated in time {\sl exactly}.  One of the advantages of this method is that there is no need to truncate the original GPE
(\ref{eq:gperot:sec5}) onto a bounded computational domain.
Again, for more details, we refer to
\cite{BaoCai313,BaoLiShen}.

Very recently, a simple and efficient numerical method has been proposed
for discretizing the GPE (\ref{eq:gperot:sec5})
\textit{via} a rotating Lagrangian coordinate \cite{Antonelli2012,BaoMTZ,Garcia2001}.
For any time $t\geq 0$, let ${ A}(t)$ be an orthogonal rotational matrix  defined as
\begin{equation}\label{Amatrix}
{ A}(t)=\left(\begin{array}{cc}
\cos(\Omega t) & \sin(\Omega t) \\
-\sin(\Omega t) & \cos(\Omega t)
 \end{array}\right),  \quad\ \  \mbox{if \ \ $d = 2$,} \qquad \quad  \ \
         \end{equation}
         and
         \begin{equation}\label{Matrixa}
{A}(t)=\left(\begin{array}{ccc}
         \cos(\Omega t) & \sin(\Omega t) & 0 \\
         -\sin(\Omega t) & \cos(\Omega t) & 0 \\
         0 & 0  & 1
         \end{array}\right), \quad\ \  \mbox{if \ \ $d = 3$.} \qquad
\end{equation}
It is easy to verify that $A^{-1}(t) =A^T(t)$ for any $t\ge0$ and ${A}(0) = { I}$, with
 $I$ the identity matrix. For any $t\ge0$, we introduce the {\it rotating Lagrangian coordinates}
$\tilde{\xbf}$ as \cite{Antonelli2012,Garcia2001}
\begin{equation}\label{transform}
\tilde{\xbf}={A}^{-1}(t) \xbf=A^T(t)\xbf \quad \Leftrightarrow \quad \xbf= {A}(t)\tilde{\xbf},
\qquad \xbf\in {\mathbb R}^d,
\end{equation}
and denote the wave function in the new coordinates as $\phi:=\phi(t,\tilde{\xbf})$
\begin{equation}\label{transform79}
\phi(t,\tilde{\xbf}):=\psi(t,\xbf)=
\psi\left(t,{A}(t)\tilde{\xbf}\right), \qquad \xbf\in {\mathbb R}^d, \quad t\geq0.
\end{equation}
In fact, here we refer the Cartesian coordinates $(t,\xbf)$ to as the {\it Eulerian coordinates}
and $(t,\tilde{\xbf})$ to as the rotating Lagrangian coordinates for any fixed $t\ge0$.
Using the chain rule, we  obtain the  following
$d$-dimensional GPE in the rotating Lagrangian coordinates without the angular momentum rotation
term \cite{Antonelli2012,BaoMTZ}:
\begin{equation}\label{GPElag321}
i\partial_t \phi(t, \tilde{\xbf}) = \left[-\frac{1}{2}\nabla^2 + W(\tilde{\xbf},t) + \beta|\phi|^2
 \right]\phi(t,\tilde{\bf x}), \qquad \tilde{\xbf}\in{\mathbb R}^d, \quad t>0,
\end{equation}
where $W(\tilde{\xbf},t) = V(A(t)\tilde{\xbf})$ for $\tilde{\xbf}\in {\mathbb R}^d$ and $t\ge0$.
The initial data (\ref{initgpe876}) can be transformed as
\begin{equation}
\label{Initial432}
\phi(t=0,\tilde{\xbf}) = \psi(t=0,\xbf)=\psi_0(\xbf):=\phi_0(\xbf)=\phi_0(\tilde{\xbf}), \qquad
\tilde{\xbf}=\xbf\in{\mathbb R}^d.
\end{equation}
Then the GPE (\ref{GPElag321}) with the initial data (\ref{Initial432})
can be directly solved  by the TSSP presented in Section 2.
After obtaining the numerical solution $\phi(t,\tilde{\xbf})$ on a bounded computational domain,
 if it is needed to recover the original wave function $\psi(t,\xbf)$ over
a set of fixed grid points in the Cartesian coordinates $\xbf$,
one can  use the standard Fourier/sine interpolation operators from the discrete numerical solution
$\phi(t,\tilde{\xbf})$ to construct an interpolation continuous function over the bounded computational domain \cite{Boyd1992,Shen}, which can be used to compute $\psi(t,\xbf)$ over
a set of fixed grid points in the Cartesian coordinates $\xbf$ for any fixed time $t\ge0$.
For more details, we refer to \cite{BaoMTZ}.

\section{Extension to coupled NLSEs/GPEs} \label{s4}
\setcounter{equation}{0}

In many applications, e.g. multi-components BEC \cite{Bao,Bao3,BaoCai313,PitaevskiiStringari} and/or
interaction of laser beams \cite{ADP,AS,BaoZheng,N}, coupled NLSEs/GPEs have been used
for modeling different problems. For simplicity of notations,
here we only consider coupled NLSEs/GPEs with two equations and cubic nonlinearity
for two-components BEC and/or interaction of two laser beams.
Extensions to coupled NLSEs/GPEs with more than two equations
are straightforward. Consider \cite{Aft,Bao,Bao01,BaoCai313,BaoP,PitaevskiiStringari,ZhangBL}
\begin{equation}\label{eq:cgpe1:sec9}
\begin{split}
&i\partial_t \psi_1(t,\xbf)=\left[-\frac 12\nabla^2
+V_1(\xbf)+\beta_{11}|\psi_1|^2+\beta_{12}|\psi_2|^2\right]\psi_1+\lambda
\psi_2, \qquad \xbf\in{\mathbb R}^d,  \quad t>0,\\
&i\partial _t \psi_2(t,\xbf)=\left[-\frac 12\nabla^2
+V_2(\xbf)+\beta_{21}|\psi_1|^2+\beta_{22}|\psi_2|^2\right]\psi_2+\lambda
\psi_1,\qquad \xbf\in {\mathbb R}^d,  \quad t>0,\end{split}
\end{equation}
with initial data
\begin{equation}
\psi_1(t=0,\xbf)=\psi_1^{(0)}(\xbf), \qquad \psi_2(t=0,\xbf)=\psi_2^{(0)}(\xbf), \qquad \xbf\in {\mathbb R}^d.
\end{equation}
Here $(\psi_1,\psi_2):=(\psi_1(\xbf,t),\psi_2(\xbf,t))$ is the dimensionless
complex-valued macroscopic wave function, $V_1(\xbf)$  and $V_2(\xbf)$ are two given dimensionless
real-valued external potentials, $\beta_{11}$, $\beta_{12}=\beta_{21}$ and $\beta_{22}$ are
given dimensionless real constants describing the interaction strength,
 and $\lambda$ is a given dimensionless
real constant describing internal atomic Josephson junction in a two-components BEC \cite{Bao,Bao01,BaoCai313,PitaevskiiStringari,Wang,WangXu,ZhangBL}.
This coupled NLSEs/GPEs conserves the {\sl total mass} as \cite{Bao,Bao01,BaoCai313,PitaevskiiStringari,ZhangBL}
 \begin{equation}\label{eq:mass1:sec9}
N(t)=\int_{{\mathbb R}^d} \left(|\psi_1(t,\xbf)|^2+|\psi_2(t,\xbf)|^2\right)d\xbf\equiv
\int_{{\mathbb R}^d} \left(|\psi_1^{(0)}(\xbf)|^2+|\psi_2^{(0)}(\xbf)|^2\right)d\xbf:=N(0),
\qquad t\ge0, \end{equation}
and the {\sl energy} as
\begin{eqnarray}
E(t)&:=&\int_{{\mathbb R}^d}\biggl[\frac
12 \left(|\nabla\psi_1|^2+
|\nabla\psi_2|^2\right)+V_1(\xbf)|\psi_1|^2+V_2(\xbf)|\psi_2|^2
+\frac 12 \beta_{11}|\psi_1|^4\nonumber\\
&&\qquad +\frac 12\beta_{22}|\psi_2|^4
+\beta_{12}|\psi_1|^2|\psi_2|^2+2\lambda\cdot\text{Re}
(\psi_1\bar{\psi}_2)\biggl]d\xbf\equiv E(0),\qquad t\ge0, \label{eq:energy:sec9}
\end{eqnarray}
where $Re(f)$ denotes the real part of the function $f$.
 In addition, if there is no internal
Josephson junction  in (\ref{eq:cgpe1:sec9}), i.e. $\lambda=0$, the mass of each
component is also conserved \cite{Bao,Bao01,BaoCai313,PitaevskiiStringari,ZhangBL}
\begin{equation}\label{eq:Njtt1:sec9} N_1(t):=
\int_{{\mathbb R}^d} |\psi_1(t,\xbf)|^2\,d\xbf\equiv \int_{{\mathbb R}^d} |\psi_1^{(0)}(\xbf)|^2\,d\xbf, \qquad
N_2(t):=\int_{{\mathbb R}^d}
|\psi_2(t,\xbf)|^2\,d\xbf\equiv \int_{{\mathbb R}^d}
|\psi_2^{(0)}(\xbf)|^2\,d\xbf, \qquad t\ge0.\end{equation}

Different efficient and accurate numerical methods have been proposed in the
 literatures for discretizing the above coupled NLSEs/GPEs \cite{Bao,Bao01,BaoCai313,Wang,WangXu,ZhangBL}.
In fact, the extension of  the numerical methods TSSP, CNFD, ReFD, SIFD and TSFD for the
NLSE/GPE (\ref{eqSchrodinger}) presented in  Section 2
is direct for the coupled NLSEs/GPEs (\ref{eq:cgpe1:sec9}).
For simplicity of notations, here we only present the TSSP method
for (\ref{eq:cgpe1:sec9}).

From time $t=t_n$ to time $t=t_{n+1}$, the coupled  NLSEs/GPEs (\ref{eq:cgpe1:sec9}) is solved in two steps.
One solves
\begin{equation}\label{eq:cgpe1:sec946}
\begin{split}
&i\partial_t \psi_1(t,\xbf)=-\frac 12\nabla^2\psi_1(t,\xbf)+\lambda
\psi_2(t,\xbf), \qquad   t>t_n,\\
&i\partial _t \psi_2(t,\xbf)=-\frac 12\nabla^2\psi_2(t,\xbf)+\lambda
\psi_1(t,\xbf),\qquad   t>t_n,\end{split}
\end{equation}
for one time step of length $\tau$, followed by solving
\begin{equation}\label{eq:cgpe1:sec948}
\begin{split}
&i\partial_t \psi_1(t,\xbf)=\left[V_1(\xbf)+\beta_{11}|\psi_1(t,\xbf)|^2+
\beta_{12}|\psi_2(t,\xbf)|^2\right]\psi_1(t,\xbf), \qquad  t>t_n,\\
&i\partial _t \psi_2(t,\xbf)=\left[V_2(\xbf)+\beta_{21}|\psi_1(t,\xbf)|^2+
\beta_{22}|\psi_2(t,\xbf)|^2\right]\psi_2(t,\xbf),\qquad  t>t_n,\end{split}
\end{equation}
for the same time step $\tau$. Similar to (\ref{NLSE1st}), Eq. (\ref{eq:cgpe1:sec946})
can be discretized in space by sine spectral method and then integrated in time {\sl exactly}
\cite{Bao,Bao3,BaoCai313,ZhangBL}. Like for (\ref{NLSE2nd}), for $t\in[t_n,t_{n+1}]$, Eq. (\ref{eq:cgpe1:sec948})
leaves $\rho_{_1}:=|\psi_1|^2$ and $\rho_{_2}:=|\psi_2|^2$ invariant \cite{Bao,BaoCai313,ZhangBL},
i.e.  $\rho_{_1}(t,\xbf)\equiv \rho_{_1}(t_n,\xbf)$ and
$\rho_{_2}(t,\xbf)\equiv \rho_{_2}(t_n,\xbf)$ for $t_n\le t\le t_{n+1}$ for any fixed $\xbf$;
thus (\ref{eq:cgpe1:sec948}) can be integrated in time {\sl exactly} \cite{Bao,Bao3,BaoCai313,ZhangBL}.
Then, we can construct the second-order TSSP
 method for the coupled NLSEs/GPEs (\ref{eq:cgpe1:sec9}) \textit{via} the Strang
 splitting \cite{Bao,Bao01,BaoCai313,BaoZhang2,ZhangBL}. We omit the details  here for brevity.

We remark here that the above TSSP, CNFD, SIFD methods have
been extended to solve many other nonlinear dispersive partial differential
equations arising from different applications.
For details, we refer to the Schr\"{o}dinger equation with wave operator \cite{BaoCai,BaoCai3},
the Sch\"{o}dinger-Poisson system \cite{BaoJinP,BMS,SPMCompare},
the Zakharov system \cite{BaoSun,BaoSun2,Chang,Glassey,Jin967},
the Klein-Gordon-Schr\"{o}dinger equations \cite{BaoY},
and the Ginzburg-Laudan-Schr\"{o}dinger equations \cite{ZhangBD} and references therein.

\section{Numerical comparison and applications}  \label{s5}
\setcounter{equation}{0}

\subsection{Comparison of different numerical methods}

In order to compare the numerical performance and accuracy
of different numerical methods,
such as TSSP, CNFD, SIFD, ReFD and TSFD, for the NLSE/GPE (\ref{eqSchrodinger}),
we take $d=1$, $\varepsilon=1$, $V(x)\equiv 0$ and $f(\rho)=-\rho$ in (\ref{eqSchrodinger}),
and the initial data  $\psi_0$ in (\ref{init0}) as
\[\psi_0(x)= A\;\text{sech}(A(x-x_0))\,e^{i(vx+\theta_0)}, \qquad x\in {\mathbb R},\]
with $A=2$, $v=1$ and $x_0=\theta_0=0$. Then the NLSE/GPE (\ref{eqSchrodinger}) with (\ref{init0})
has the exact bright soliton solution (\ref{bright_sol}), i.e.
$\psi(t,x)=\psi_{B}(t,x)$, with $\beta=-1$, $A=2$, $v=1$ and $x_0=\theta_0=0$.
In our computation, we take the bounded computational domain as the interval $(a,b)$ with
$a=-15$ and $b=20$ and homogeneous Dirichlet boundary condition, which are
large enough so that the truncation errors can be ignored.
In order to quantify the numerical solution, we use
the  $l^{\infty}$-norm of the error between
the numerical solution $\psi_j^n$ and the exact solution
$\psi(t_n,x_j)$ as
\begin{equation}\label{l_err}
  e_{\infty}^p(t_n):=\max_{0\le j\le J}|\psi(t_n,x_j)-\psi_{j}^{n}|, \qquad
 e_{\infty}^m(t_n):=\max_{0\le j\le J} ( |\psi(t_n,x_j)|-|\psi_{j}^{n}|), \qquad
n\ge0.
\end{equation}
The functions $e_{\infty}^p$ and $e_{\infty}^m$ allow respectively to measure the phase error and the modulus error.\\

\begin{table}[htbp]
\begin{center}
\begin{tabular}{|c c|c|c|c|c|c|}
  \hline
$h$ & & $h_0$=0.5   & $h_0$/2   &$h_0$/4     &  $h_0$/8   & $h_0$/16 \\ \hline
 CNFD &$e_{\infty}^p$&  2.48     &  1.87E0        & 4.28E-1 & 1.03E-1 & 2.57E-2 \\
      &$e_{\infty}^m$ &  1.89     &  6.98E-1     & 1.46E-1 & 3.52E-2 & 8.68E-3 \\ \hline
 ReFD &$e_{\infty}^p$&  2.48     &  1.87E0        & 4.28E-1 & 1.03E-1 & 2.57E-2 \\
      &$e_{\infty}^m$ &  1.89     &  6.98E-1     & 1.46E-1 & 3.52E-2 & 8.68E-3 \\ \hline
 SIFD &$e_{\infty}^p$&  2.48     &  1.87E0        & 4.28E-1 & 1.03E-1 & 2.57E-2 \\
      &$e_{\infty}^m$ &  1.89     &  6.98E-1     & 1.46E-1 & 3.52E-2 & 8.68E-3 \\ \hline
 TSFD &$e_{\infty}^p$&  2.48     &  1.87E0        & 4.28E-1 & 1.03E-1 & 2.57E-2 \\
      &$e_{\infty}^m$ &  1.89     &  6.98E-1     & 1.46E-1 & 3.52E-2 & 8.68E-3 \\ \hline
 TSSP &$e_{\infty}^p$ &  1.485     &  3.81E-4     & 8.63E-9 & $<$1E-9 & $<$1E-9 \\
      &$e_{\infty}^m$ &  1.408     &  2.45E-4     & 4.49E-9 & $<$1E-9 & $<$1E-9 \\
\hline
\end{tabular}
\end{center}
 \caption{Spatial error analysis on errors $e_{\infty}^{p,m}(t=5)$  of
different numerical methods for the NLSE/GPE (\ref{eqSchrodinger}) in 1D under different mesh
sizes $h$.  }
\label{table_comp_bright_dir}
\end{table}

To test the spatial discretization errors of the different numerical methods,
we fix $\tau=10^{-5}$ such that the time discretization errors are negligible.
Table \ref{table_comp_bright_dir} shows
the spatial  errors $e_{\infty}^{p,m}(t=5)$ for different
 numerical methods under different mesh sizes $h$.
Similarly,  to compare the temporal discretization errors of different numerical methods,
we take $h=3.5 \times 10^{-3}$ to get very small  spatial discretization errors.
Table \ref{table_comp_bright_tt} shows
the temporal  errors $e_{\infty}^{p,m}(t=5)$ for different
 numerical methods under different time steps $\tau$ \cite{BaoTangXu}.


\begin{table}[htbp]
\begin{center}
\begin{tabular}{|c c|c|c|c|c|c|}
  \hline
$\tau$& & $\tau_0$=0.1  & $\tau_0$/2 & $\tau_0$/4 & $\tau_0$/8 & $\tau_0$/16 \\ \hline
 CNFD &$e_{\infty}^p$ & 2.62E-1 & 6.65E-2 & 1.64E-2 & 3.88E-3 & 7.28E-4 \\
      &$e_{\infty}^m$ & 1.08E-2 & 2.87E-3 & 6.70E-4 & 1.16E-4 & 6.64E-5 \\ \hline
 ReFD &$e_{\infty}^p$ & 3.11E-1 & 2.68E-2 & 1.97E-2 & 5.07E-3 & 1.47E-3 \\
      &$e_{\infty}^m$ & 2.25E-1 & 5.37E-2 & 1.34E-2 & 3.37E-3 & 9.25E-4 \\ \hline
 SIFD &$e_{\infty}^p$ & 2.96E0  & 6.46E-1 & 1.91E-1 & 6.62E-2 & 2.61E-2 \\
      &$e_{\infty}^m$ & 2.57E-1 & 1.13E-1 & 5.34E-2 & 2.58E-2 & 1.26E-2 \\ \hline
 TSFD &$e_{\infty}^p$ & 8.51E-1 & 2.00E-1 & 4.97E-2 & 1.25E-2 & 3.37E-3 \\
      &$e_{\infty}^m$ & 4.10E-1 & 9.03E-2 & 2.20E-2 & 5.49E-3 & 1.44E-3 \\ \hline
 TSSP &$e_{\infty}^p$ & 5.17E-1 & 1.40E-1 & 3.57E-2 & 8.98E-3 & 2.25E-3 \\
      &$e_{\infty}^m$ & 4.98E-2 & 1.64E-2 & 4.21E-3 & 1.06E-3 & 2.65E-4 \\
\hline
\end{tabular}
\end{center}
 \caption{Temporal error analysis on errors $e_{\infty}^{p,m}(t=5)$  of
different numerical methods for the NLSE/GPE (\ref{eqSchrodinger}) in 1D under different time
steps $\tau$.  }
\label{table_comp_bright_tt}
\end{table}


From Tabs. \ref{table_comp_bright_dir} \& \ref{table_comp_bright_tt} as well as additional
numerical results not shown here for brevity, it is clearly demonstrated that
TSSP is spectral order accurate in space and second-order accurate in time while
 CNFD, ReFD, SIFD and TSFD are second-order accurate in both space and time.
For more numerical comparisons, we refer to \cite{BaoCai313,BaoTangXu,Chang1} and references therein.

\subsection{Applications}

\begin{figure}[htb!]
\centerline{\psfig{figure=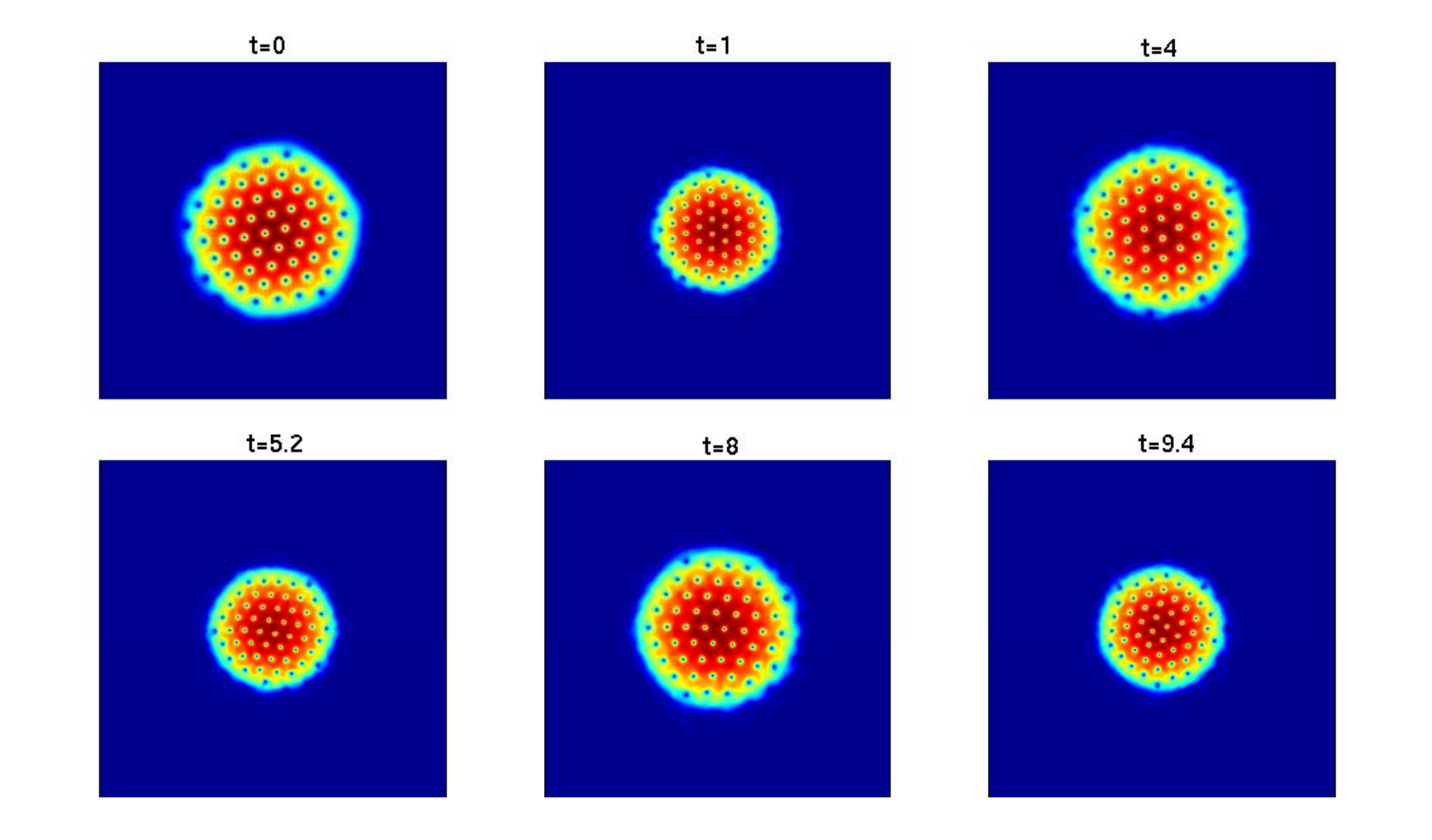,height=9cm,width=15.5cm,angle=0}}
\centerline{\psfig{figure=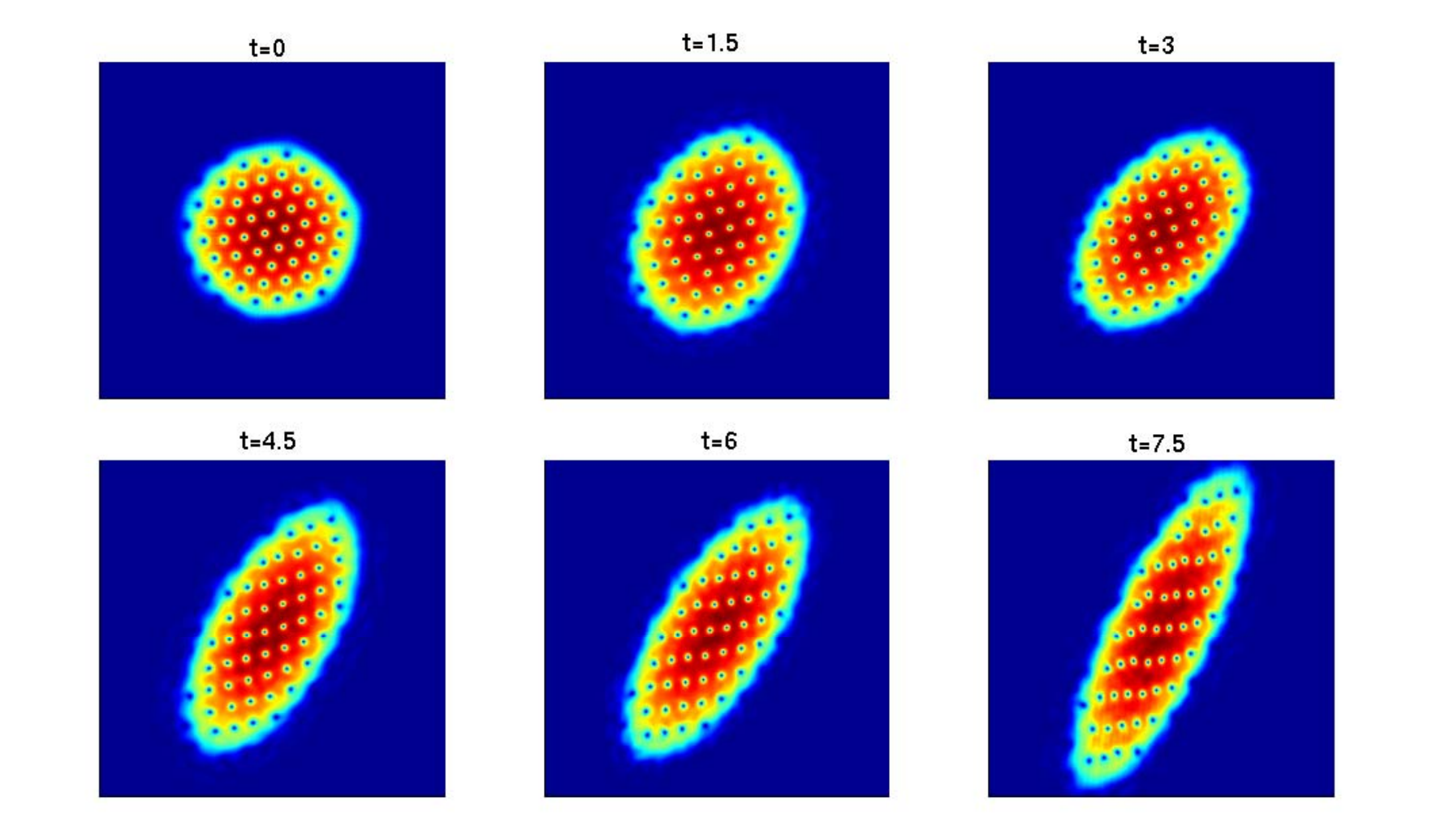,height=9cm,width=15.5cm,angle=0}}
\caption{Contour plots of the density function $\rho:=|\psi(t,\xbf)|^2$ for the
dynamics of a quantized vortex lattice in a rotating
BEC for case I  (top two rows)
and case II  (bottom two rows).}\label{F5}
\end{figure}

In order to show numerical results for problems coming from applications,
we consider the dynamics of the NLSE/GPE (\ref{eq:gperot:sec5}) with a rotation
term in 2D  starting from a quantized vortex lattice for a rotating BEC \cite{BaoCai313,BaoCai2,BaoLiShen,BaoMTZ,BaoWang},
i.e. we take $d=2$,  $\beta = 1000$ and $\Omega = 0.9$ in (\ref{eq:gperot:sec5}). The initial datum
in (\ref{initgpe876}) is chosen as a
stationary vortex lattice which is computed numerically by using
the method in \cite{BaoCai313,BaoWangP} with the above parameters and
the harmonic potential $V(x,y)=\frac{1}{2}(x^2+y^2)$. Then, the dynamics of the vortex lattice is
studied numerically by perturbing the
harmonic potential from $V(x,y)=\frac{1}{2}(x^2+y^2)$
to $V(x,y)=\frac{1}{2}(\gamma_x^2 x^2+\gamma_y^2 y^2)$
with (i) case I: $\gamma_x=\gamma_y=2$, and (ii) case II: $\gamma_x=1.1$ and $\gamma_y=0.9$,
respectively, at time $t=0$.
In the numerical simulation, we choose the bounded computational domain
as $[-32, 32]^2$ with homogeneous Dirichlet boundary condition,
mesh size $h =1/16$ and time step $\tau =10^{-4}$. Figure \ref{F5}  shows the contour plots of
the density function $\rho(t,\xbf):=|\psi(t,\xbf)|^2$ displayed on $[-17, 17]^2$
at different times for cases I and II.

\section{Conclusion and perspectives} \label{s6}
\setcounter{equation}{0}

Due to its massive applications in many different areas, the research on
numerical methods and simulation as well as applications related to the
nonlinear Schr\"{o}dinger/Gross-Pitaevskii equations (NLSE/GPE)
has been started several decades ago. Up to now, rich and extensive research results
have been obtained in developing and analyzing efficient and accurate
numerical methods for the NLSE/GPE and in applying them for simulating problems
arising from many different areas, such as Bose-Einstein condensation (BEC),
nonlinear optics, superfluids, etc. Nowadays, numerical simulation has
become a very important tool in theoretical and computational physics
as well as computational and applied mathematics for solving problems
related to NLSE/GPE and it has been used to predict and guide new experiments
due to the advances in numerical methods and their analysis as well
as powerful and/or parallel computers.
Of course, due to its dispersive nature of the NLSE/GPE, when one is doing
numerical simulation, he/she should choose the numerical method,
mesh size and time step as well as the computational domain properly and carefully
so that the numerical results reflect ``correct'' physical phenomena.

The research in this area
is still very active and highly demanded due to the latest
experimental and/or technological advances in BEC,
nonlinear optics, graphene, semiconductors, topological insulators,
materials simulation and design, etc.
It becomes more and more interdisciplinary involving
theoretical, computational and experimental physicists and
computational and applied mathematicians as well as computational scientists.
Of course, there are still many important issues to be addressed.
For example, extension and designing as well as analyzing
new numerical methods for highly oscillatory nonlinear dispersive
partial differential equations,
especially coupled NLSE/GPE with other equations such as the
Davey-Stewartson system \cite{SulemSulem},
Kadomstev-Petiashvili equations \cite{SulemSulem}, coupled NLSE/GPE with quantum
Boltzmann equation for BEC at finite temperature \cite{BaoCai313,Zaremba}, etc., are always welcome
and highly demanded. Another issue is to design efficient and accurate numerical methods
and apply them for studying  numerically NLSE/GPE with random potential \cite{And,Dubi,Lye,MinLRB}
or stochastic NLSE/GPE \cite{Barton-SmithDebisscheDiMenza, BouardDebussche, DeBouardDebusscheDiMenza,   DebusscheDiMenza,  Marty}
with applications and NLSE/GPE in higher dimensions, i.e. $d>3$
for many-body problems in quantum chemistry and materials simulation and design.
Here memory and computational cost might be extremely high and thus parallel computing
and/or sparse grids as well as spatial and temporal adaptivity
are very useful and essential. Last but not the least, efficient implementation
of numerical methods and portable and readable programming are very important
from the application point of view. Although there are plenty of research codes available
in the community \cite{CaliariRainer,Caplan,MurA,VVBMA},
a software package or tool box (with parallel implementation in 2D and 3D)
is still very useful for applications. In summary, in order to solve
challenging scientific and engineering problems
and/or guiding and predicting new experiments related to NLSE/GPE,
the interaction and close collaboration between computational and applied mathematicians
and theoretical and experimental physicists  and chemists as well
as computational and applied scientists becomes more and more essential
in this research area.

\bigskip

\noindent {\bf Acknowledgements.}

This work was partially supported by the French ANR
grants MicroWave NT09\_460489 (''Programme Blanc'' call) and
ANR-12-MONU-0007-02 BECASIM (''Mod\`eles Num\'eriques''  call) (X. Antoine and C. Besse),
and  by the Singapore A*STAR SERC  PSF-Grant	1321202067 (W. Bao).

\end{document}